\documentclass[11pt]{article}
\usepackage{amsmath,amssymb,amsthm,bm,graphicx,color,epsfig}

\theoremstyle{plain}
\newtheorem{theorem}{Theorem}
\newtheorem{lemma}{Lemma}

\theoremstyle{definition}

\newcommand{\R}{\mathbb{R}}
\newcommand{\Z}{\mathbb{Z}}

\newcommand{\Grad}{\nabla}

\newcommand{\dif}{\mathrm{d}}
\newcommand{\pd}{\partial}
\renewcommand{\d}{\dif}
\newcommand{\eps}{\varepsilon}

\newcommand{\hxwa}{\hat{\xi}_\ell |\xi|}
\newcommand{\hxw}{\hat{\xi}_\ell}
\newcommand{\tw}{\theta_\ell}

\numberwithin{equation}{section}

\title{Fast Computation of Fourier Integral Operators}

\author{Emmanuel Cand\`{e}s$^\dagger$, Laurent Demanet$^\sharp$ and Lexing Ying$^\ast$\\
  [2mm]
  $\dagger$ Applied and Computational Mathematics, Caltech, Pasadena, CA 91125\\
  $\sharp$ Department of Mathematics, Stanford University, Stanford, CA 94305\\
  $\ast$ Department of Mathematics, University of Texas, Austin, TX 78712
}
\date{\today}

\begin{document}
\maketitle

\begin{abstract}
  We introduce a general purpose algorithm for rapidly computing
  certain types of oscillatory integrals which frequently arise in
  problems connected to wave propagation and general hyperbolic
  equations. The problem is to evaluate numerically a so-called
  Fourier integral operator (FIO) of the form $\int e^{2\pi i
    \Phi(x,\xi)} a(x,\xi) \, \hat{f}(\xi) \mathrm{d}\xi$ at points
  given on a Cartesian grid. Here, $\xi$ is a frequency variable,
  $\hat f(\xi)$ is the Fourier transform of the input $f$, $a(x,\xi)$
  is an amplitude and $\Phi(x,\xi)$ is a phase function, which is
  typically as large as $|\xi|$; hence the integral is highly
  oscillatory at high frequencies. Because an FIO is a dense matrix, a
  naive matrix vector product with an input given on a Cartesian grid
  of size $N$ by $N$ would require $O(N^4)$ operations.

  This paper develops a new numerical algorithm which requires
  $O(N^{2.5} \log N)$ operations, and as low as $O(\sqrt{N})$ in
  storage space. It operates by localizing the integral over polar
  wedges with small angular aperture in the frequency plane. On each
  wedge, the algorithm factorizes the kernel $e^{2 \pi i \Phi(x,\xi)}
  a(x,\xi)$ into two components: 1) a diffeomorphism which is handled
  by means of a nonuniform FFT and 2) a residual factor which is
  handled by numerical separation of the spatial and frequency
  variables. The key to the complexity and accuracy estimates is that
  the separation rank of the residual kernel is \emph{provably
    independent of the problem size}. Several numerical examples
  demonstrate the efficiency and accuracy of the proposed methodology.
  We also discuss the potential of our ideas for various applications
  such as reflection seismology.
\end{abstract}


{\bf Acknowledgments.} E.~C. is partially supported by an NSF grant
CCF-0515362 and a DOE grant DE-FG03-02ER25529.  L.~D. and L.~Y. are
supported by the same NSF and DOE grants. We are thankful to William
Symes for stimulating discussions about Kirchhoff migration and
related topics.


\section{Introduction}
\label{sec:intro}

This paper introduces a general-purpose algorithm to compute the
action of linear operators which are frequently encountered in
analysis and scientific computing. These operators take the form
\begin{equation}
\label{eq:FIO}
  (Lf)(x) = \int_{\R^d} a(x,\xi) e^{2 \pi i \Phi(x,\xi)} \hat{f}(\xi) 
\, \d\xi, 
\end{equation}
where $\Phi(x,\xi)$ is a smooth phase function obeying the homogeneity
relation $\Phi(x,\lambda \xi) = \lambda \Phi(x,\xi)$ for $\lambda$
positive, and $a(x,\xi)$ is a smooth amplitude term. As is standard,
$\hat f$ is the Fourier transform of $f$ defined by
\begin{equation}
  \label{eq:fourier}
  \hat f(\xi) = \int_{\R^d} f(x) e^{-2 \pi i x \xi} \, \d x.  
\end{equation}
With the proper regularity assumptions on the phase and amplitude to
be detailed later, \eqref{eq:FIO} defines a class of oscillatory
integrals known as Fourier integral operators (FIOs). FIOs are the
subject of considerable study for many of the operators encountered in
physics and other fields are of this form. For instance, most
differential and pseudodifferential operators are FIOs.  Convolutions
and multiplications by smooth functions are FIOs. Some ``principal
value'' integrals are FIOs. And the list goes on.

An especially important example of FIO is the solution operator to the
free-space wave equation in $\R^d$, $d > 1$,
\begin{equation}\label{eq:wave}
\frac{\pd^2 u}{\pd t^2}(x,t) = c^2 \Delta u(x,t), 
\end{equation}
with initial conditions $u(x,0) = u_0(x)$ and $\frac{\pd u}{\pd
  t}(x,0) = 0$, say. Everyone knows that for constant speeds, the
Fourier transform decouples the different frequency components of
$u$. Each Fourier component obeys an ordinary differential equation
which can be solved explicitly. The solution $u(x,t)$ is the
superposition of these Fourier modes and is given by
\begin{equation}\label{eq:wave_solution}
u(x,t) = \frac{1}{2}\left(
\int e^{2\pi i (x\cdot\xi + c |\xi| t)} \hat{u_0}(\xi) \d\xi +
\int e^{2\pi i (x\cdot\xi - c |\xi| t)} \hat{u_0}(\xi) \d\xi
\right).
\end{equation}
The connection is now clear: the solution operator is the sum of two
Fourier integral operators with phase functions
\[
\Phi_{\pm}(x,\xi) = x\cdot\xi \pm c |\xi| t.
\]
For variable but reasonably smooth sound speeds $c(x)$, the solution
operator is for small times a sum of two FIOs with more
complicated phases and amplitudes. In particular, the phase can be
constructed from the optical traveltime in a medium with index of
refraction $1/c(x)$, see \cite{candes-2005-crwpos} for details.

In short, it is useful to think of FIOs as proxies for the solution
operator to large classes of hyperbolic differential equations.

\subsection{FIO computations}

Numerical simulation of free wave propagation with constant sound
speed is straightforward. As long as the solution $u(x,t)$ is
sufficiently well localized both in space and frequency, it can be
computed accurately and rapidly by applying the sequence of steps
below.
\begin{enumerate}
\item Compute the Fast Fourier Transform (FFT) of $u_0$. 
\item Multiply the result by $e^{\pm 2 \pi i c |\xi| t}$, and sum as in
  (\ref{eq:wave_solution}). 
\item Compute the inverse FFT.
\end{enumerate}
Of course, this only works in the very special case where the
amplitude $a$ is independent of $x$, and where the phase is of the
form $x \cdot \xi$ plus a function of $\xi$ alone. Expressed
differently, this works when the FIO is shift-invariant so that it is
diagonal in the Fourier basis.  Note that there is in general no
formula for the eigenfunctions when $\Phi$ or $a$ depend on $x$.
Computing these eigenfunctions on the fly is out of the question when
the objective is merely to compute the action of the operator. (Note
that even if the spectral decomposition of the operator were
available, it is not clear how one would use it to speed up
computations.)

The object of this paper is to find an algorithm that is considerably
faster than evaluating \eqref{eq:FIO} by direct quadratures, and is
yet suited to handle large classes of phases and amplitudes. Most of
the existing fast summation techniques rely on either the
non-oscillatory behavior (such as wavelet based techniques
\cite{beylkin-1991-fwtna}) or the existence of a low rank
approximation (fast multipole methods \cite{greengard-1987-FAPS},
hierarchical matrices \cite{hackbusch-1999-smabhm}, pseudodifferential
separation \cite{bao-1996-cpdo}). The difficulty here is that the
kernel $e^{2\pi i \Phi(x,\xi)}$ is highly oscillatory and does not
have a low rank separated approximation. Therefore, all the modern
techniques are not directly applicable.

The main claim of this paper, however, is that there is a way to
decompose the operator into a sum of components for which the
oscillations are well-understood and low-rank representations are
available. In addition, the number of such components is reasonably
small which paves the way to faster algorithms. Before expanding on
this idea, we first explain the discretization of the operator
\eqref{eq:FIO}.

\subsection{Discretization}

For simplicity, we restrict our attention in this paper to the two
dimensional case $d=2$. The situation in which $d \ge 3$ is exactly
the same.

Just as the discrete Fourier transform is the digital analogue of the
continuous Fourier transform, one can also introduce discrete Fourier
integral operators. Given a function $f$ defined on a Cartesian grid
$X = \{x=(\frac{n_1}{N}, \frac{n_2}{N}), 0\le n_1, n_2 <N
\;\mbox{and}\;n_1,n_2 \in \Z \}$, we simply define the discrete
Fourier integral operator by
\begin{equation}
  (Lf)(x) := \frac{1}{N} \sum_{\xi \in \Omega} a(x,\xi) e^{2\pi i \Phi(x,\xi)} \hat{f}(\xi)
  \label{eq:dfio}
\end{equation}
for every $x\in X$.  (We are sorry for overloading the symbol $L$ to
denote both the discrete and continuous object but there will be no
confusion in the sequel.) The summation above is taken over all
$\Omega = \{\xi = (n_1,n_2), -\frac{N}{2} \le n_1,n_2
<\frac{N}{2}\;\mbox{and}\;n_1,n_2 \in \Z \}$ and throughout this
paper, we will assume that $N$ is an even integer.  Here and below,
$\hat{f}$ is the discrete Fourier transform (DFT) of $f$ and is
defined as
\begin{equation}
  \label{eq:dft}
  \hat{f}(\xi) = \frac{1}{N} \sum_{x \in X} e^{-2 \pi i x\cdot\xi} f(x).
\end{equation}
The normalizing constant $\frac{1}{N}$ in \eqref{eq:dfio}
(resp.~\eqref{eq:dft}) ensures that $L$ (resp.~the DFT) is a discrete
isometry in the case where $\Phi(x,\xi) = x \cdot \xi$.

The formula \eqref{eq:dfio} turns out to be an accurate discretization
of (\ref{eq:FIO}) as soon as $f$ obeys standard localization estimates
both in space and frequency. A justification of this fact would
however go beyond the scope of this paper, and is omitted.  In the
remainder of the paper, we will take (\ref{eq:dfio}) as the quantity
we wish to compute once we are given a phase and an amplitude function.

The parameter $N$ measures the size and difficulty of the
computational problem. In a nutshell, it corresponds to the number of
points which are needed in each direction to accurately sample the
continuous object $f(x)$. This is the reason why $N$ will be a central
quantity throughout the rest of paper.

As mentioned earlier, the straightforward method for computing
\eqref{eq:dfio} simply evaluates the summation independently for each
$x$. Since each sum takes $O(N^2)$ operations and there are $N^2$ grid
points in $X$, this strategy requires $O(N^4)$ operations. When $N$ is
moderately large, this can be prohibitive.  This paper describes a
novel algorithm which computes all the values of $Lf(x)$ for $x \in X$
with high accuracy in $O(N^{2.5} \log N)$ operations. The only
requirement is that the amplitude and the phase obey mild smoothness
conditions, which are in fact standard.

\subsection{Separation within angular wedges}

\begin{figure}[htb]
  \begin{center}
        \includegraphics[height=2in]{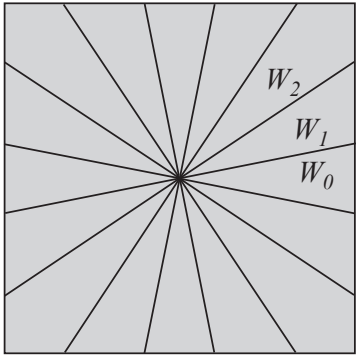}
  \end{center}
  \caption{The frequency domain is partitioned into $\sqrt{N}$
    equiangular wedges.}
  \label{fig:wedge}
\end{figure}

This section outlines the main idea of the paper. Let $\arg \xi$ be
the angle between $\xi$ and the horizontal vector $(1,0)$, and partition
the frequency domain into a family of angular wedges $\{W_\ell\}$
defined by
\[
W_\ell = \{\xi : {(2\ell-1)\pi}/{\sqrt{N}} \le \arg \xi <
{(2\ell+1)\pi}/{\sqrt{N}} \}
\]
for $0 \le \ell < \sqrt{N}$ (assume $\sqrt{N}$ is an integer). An
important property of these wedges is that each $W_\ell$ satisfies the
parabolic relationship
\begin{equation}
  length \simeq width^2,
  \label{eq:parabolic}
\end{equation}
up to multiplicative constants independent of $N$. There are
$O(\sqrt{N})$ such wedges, as illustrated in Figure
\ref{fig:wedge}. 

For each wedge $W_\ell$, we let $\chi_\ell$ be the indicator function
of $W_\ell$. Similarly, we denote by $\hat \xi_\ell$ the unit vector
pointing to the center direction of $W_\ell$
\[
\hat\xi_\ell = \left( \cos \frac{2\ell\pi}{\sqrt{N}}, 
  \sin \frac{2\ell\pi}{\sqrt{N}} \right).
\]
It follows from the identify $\sum_{\ell} \chi_\ell(\xi) = 1$ that one
can decompose the operator $L$ as $\sum_\ell L_\ell$, where
\[
(L_\ell f)(x) = \frac{1}{N} 
\sum_\xi a(x,\xi) e^{2\pi i \Phi(x,\xi)} \chi_\ell(\xi) \hat{f}(\xi).
\]
Within each wedge $W_\ell$, we can perform a Taylor expansion of
$\Phi(x,\xi)$ in the second variable, around the point $\hxwa$. There
is a point $\xi^\star$ which belongs to the line segment $[\hxwa,
\xi]$ such that
\[
\Phi(x,\xi) = \Phi(x,\hxwa) + \Grad_\xi \Phi(x,\hxwa) \cdot
(\xi-\hxwa) + \frac{1}{2} (\xi^\star-\hxwa)^T \Grad_{\xi\xi}
\Phi(x,\xi^\star) (\xi^\star-\hxwa).
\]
By homogeneity of the phase ($\Phi(x,\lambda \xi) = \lambda
\Phi(x,\xi)$ for $\lambda>0$), it holds that $\Phi(x,\xi) = \xi \cdot
\Grad_\xi \Phi(x,\xi)$ and $\Grad_\xi \Phi(x,\xi) = \Grad_\xi
\Phi(x,\hat{\xi})$. The first and third terms in the above expression
cancel and thus
\[
\Phi(x,\xi) = \Grad_\xi \Phi(x,\hxw)\cdot\xi + \frac{1}{2}
(\xi^\star-\hxwa)^T \Grad_{\xi\xi} \Phi(x,\xi^\star)
(\xi^\star-\hxwa).
\]
The first term $\Grad_\xi \Phi(x,\hxw)\cdot\xi$, which is linear in
$\xi$, is called {\em the linearized phase} and poses no problem as we
will see later on. The rest, denoted as $\Phi_\ell(x,\xi) =
\Phi(x,\xi)-\Grad_\xi \Phi(x,\hxw)\cdot \xi$ and called {\em the
  residual phase}, is of order $O(1)$ for $\xi\in W_\ell$,
independently of $N$. This follows from
\[
\Grad_{\xi\xi} \Phi(x,\xi^\star) = O(|\xi^\star|^{-1}) =
O(|\xi|^{-1}),
\]
since $\Phi(x,\xi)$ is homogeneous of degree $1$ in $\xi$, together
with
\[
|\xi^\star -\hxwa|^2 \le |\xi-\hxwa|^2 = O(|\xi|^2/N) = O(|\xi|)
\]
for all $|\xi| \le N$, which uses the fact that the shape of $W_\ell$
obeys the parabolic relationship \eqref{eq:parabolic}.

Because the residual phase $\Phi_\ell(x,\xi)$ is of order $O(1)$
independently of $N$, we say that the function $e^{2\pi i \Phi_\ell
  (x,\xi)}$ is \emph{nonoscillatory}. Under mild assumptions, this
observation guarantees the existence of a low rank separated
representation which decouples the variables $x$ and $\xi$ and
approximates the complex exponential very well.  Define the
$\epsilon$-separation rank of a function $f(x,y)$ of two variables as
the smallest integer $r_\epsilon$ for which there exists $c_n(x)$,
$d_n(y)$ such that
\[
|f(x,y) - \sum_{n = 0}^{r_\epsilon - 1} c_n(x) d_n(y)| \leq \epsilon.
\]
Then we prove the following theorem in Section \ref{sec:analy}.

{\bf Theorem.} For all $0 < \epsilon \leq 1$, there exist $N^* > 0$
and $C > 0$ such that for all $N \geq N^*$, the $\epsilon$-separation
rank of $e^{2 \pi i \Phi_\ell(x,\xi)}$ for $x \in [0,1]^2$ and $\xi
\in W_\ell$ obeys
\begin{equation}\label{eq:rank1_intro}
r_\epsilon \leq \log_2(C \epsilon^{-1}).
\end{equation}
In Section \ref{sec:analy} we make explicit the values of the
constants $N^*$ and $C$ by relating them, among other things, to the
smoothness of $\Phi$ and the angular span of $W_\ell$. We will also
provide results in the case where $N \leq N^*$, and explain why the
separation rank for the amplitude is also under control.

The point of the theorem is that the bound on the $\epsilon$-rank does
not grow as a function of $N$---in fact, the threshold condition on
$N$ indicates that the $\epsilon$-rank \emph{decays} as $N$ grows. The
logarithmic dependence on $\epsilon$ is the signature of what is
usually called spectral accuracy.

Note that the decomposition into frequency wedges obeying the
parabolic scaling has a long history in mathematics. A multiscale
version of this partitioning, the second dyadic decomposition, was
introduced by Fefferman in 1973 for the study of Bochner-Riesz
multipliers \cite{fefferman-1973-nssm}, and used by Seeger, Sogge and
Stein in 1991 to prove a sharp $L^p$-boundedness result for FIO
\cite{seeger-1991-rpfio}. More recently, it also served as the basis
for the construction of curvelets, with applications to sparsity of
FIOs and related results for wave equations \cite{smith-1998-hsfio,
  candes-2003-cfio, candes-2005-crwpos}.

\subsection{Outline of the algorithm}
The low-rank separated representation provided by the theorem above
offers us a way to compute \eqref{eq:dfio} efficiently with high
accuracy. Each term in the decomposition $L f = \sum_\ell L_\ell f$
can be further simplified as follows:
\begin{eqnarray}
  (L_\ell f)(x) &=&  \frac{1}{N} \sum_{\xi} a(x,\xi) e^{2\pi i \Phi(x, \xi)} 
  \chi_\ell(\xi) \hat{f}(\xi) 
  \nonumber\\
  &=&  \frac{1}{N} \sum_\xi e^{2\pi i \Grad_\xi \Phi(x,\hxw)\cdot\xi} 
  \; a(x,\xi) e^{2\pi i \Phi_\ell(x,\xi)} \; \chi_\ell(\xi) \hat{f}(\xi) 
  \nonumber\\
  &=& \frac{1}{N} \sum_\xi e^{2\pi i \Grad_\xi \Phi(x,\hxw)\cdot\xi} \; 
  \sum_{t = 1}^\infty \gamma^x_{\ell t}(x) \gamma^\xi_{\ell t}(\xi) \; 
  \chi_\ell(\xi) \hat{f}(\xi) \nonumber\\
  &=& \frac{1}{N} \sum_{t = 1}^\infty \gamma^x_{\ell t}(x) 
  \sum_\xi e^{2\pi i \Grad_\xi \Phi(x,\hxw)\cdot\xi} 
  \left[  \gamma^\xi_{\ell t}(\xi) \chi_\ell(\xi) \hat{f}(\xi) \right]. 
  \label{eq:dfioalgo}
\end{eqnarray}
Our analysis guarantees that the sum over $t$ can be truncated to a
fixed, hopefully small number of terms without significant loss of
precision.

In order to carry out the final summation over $t$, we first need to
construct the functions $\gamma^x_{\ell t}(x)$ and $\gamma^\xi_{\ell
  t}(\xi)$.  Sections \ref{sec:algo-dete} and \ref{sec:algo-rand}
discuss two different methods to find these functions. In Section
\ref{sec:algo-dete} we present an elementary deterministic approach,
while in Section \ref{sec:algo-rand} we present a randomized approach
that offers better efficiency both timewise and storagewise. Assuming
that $\gamma^x_{\ell t}(x)$ and $\gamma^\xi_{\ell t}(\xi)$ are
available for all values of $\ell$ and $t$, the computation of $(Lf)$
for a given $f$ consists of the following 4 steps:
\begin{enumerate}
\item Fourier transform $f$ by means of the FFT to get $\hat{f}$. 
\item Choose a bound $q$ greater than the $\eps$-rank $r_\eps$. For
  each $\ell$ and $t \leq q$, form $\hat{f}_{\ell t}(\xi) :=
  \gamma^\xi_{\ell t}(\xi) \chi_\ell(\xi) \hat{f}(\xi)$. 
\item For each $\ell$ and $t \leq q$, compute $g_{\ell t}(x) :=
  \sum_\xi e^{2\pi i \Grad_\xi \Phi(x,\hxw)\cdot\xi} \hat{f}_{\ell
    t}(\xi)$ by means of a nonuniform FFT algorithm. 
\item Compute $(Lf)(x) \approx \frac{1}{N} \sum_\ell \sum_{t = 1}^{q}
  \gamma^x_{\ell t}(x) g_{\ell t}(x)$. 
\end{enumerate}
The only step that require further discussion is the computation of
$g_{\ell,t}$. We defer the details to Section \ref{sec:algo-eval}.

It is instructive to understand why linearizing the phase is so
important. If we disregard the error introduced by the discretization
in $\xi$, we observe that $g_{\ell,t}(x)$ is simply
\[
g_{\ell,t}(x) = f_{\ell,t}(\Grad_\xi \Phi(x,\hxw)).
\]
The interpretation of an oscillatory integral in the Fourier domain as
a diffeomorphism is only possible when the phase is linear in
$\xi$. For each $\ell$ and $t$, the computation of $g_{\ell,t}$ which
is an interpolation problem, is therefore much simpler problem than
applying the original operator. Admittedly, diffeomorphisms do not
provide accurate approximations to FIOs over angular wedges, but the
content of our analysis in Section \ref{sec:analy} shows that the
computational budget to make up for the residual is safely under
control.

\subsection{Significance}

Applying nontrivial FIOs repeatedly is a daunting task that has proved
to be the computational bottleneck in various inverse problems. There
is serious scientific as well as industrial interest in speeding up
FIO computations, and accordingly a lot of resources have been
invested over the past decades in engineering better codes.

We believe that the ideas introduced in this paper provide new
directions.  To explain and illustrate this contrast, let us consider
an example from the field of reflection seismology: Kirchhoff
migration. The problem is to produce an image of the discontinuities
in the Earth's upper crust from seismograms, i.e., wave measurements
$f(t,x_r)$ parameterized by time $t$ and receiver coordinate
$x_r$. Glossing over the details, the core of Kirchhoff migration
consists in integrating several different functions $f(t,x_r)$ over a
fixed set of curves, parameterized as the level lines of some
traveltime function $\tau(x,x_r) + \tau(x,x_s)$:
\[
g(x) = \int \delta(t - \tau(x,x_r) - \tau(x,x_s)) f(t,x_r) \, \d t \, \d x_r,
\]
where $x_s$ is for us a fixed parameter (the source coordinate). We do
not expect the reader unfamiliar with seismic imaging to understand
all the physics underlying this equation. Anyone interested in details
may want to consult \cite{symes-1998-mrs}, for example. This collection
of integrals is called a generalized Radon transform (GRT), or in the
field of image processing, a Hough transform. (For convenience, the
Appendix explains why integration along ellipses---a simple GRT---is a
sum of two FIOs.) A useful notation for Kirchhoff migration is $g(x) =
(F^* f)(x)$, where $F^*$ is called the imaging operator.

The standard algorithm for applying the imaging operator is a simple
quadrature of $f(t,x_r)$, interpolated and integrated along each curve
$t = \tau(x,x_r) + \tau(x,x_s)$ (parameterized by $x$.) If the data
$f(t,x_r)$ oscillates at a wavelength comparable to the grid spacing
$1/N$, then an accurate quadrature on a smooth curve requires $O(N)$
points. Since $x$ takes on $O(N^2)$ values, the curve integration
results in a total complexity of $O(N^3)$ for applying the imaging
operator (which is of course better than the $O(N^4)$ complexity of
the naive summation.)

In reality, the true $F^*$ is only approximated by a GRT. The
derivation of the expression for $F^*$ from the wave equation reveals
that if the geometry of the optical rays is not too complex, $F^*$ is
in fact closer to an FIO than a GRT \cite{symes-1998-mrs}. This is
akin to the observation that the retarded propagator of the wave
equation in 2D is not a distribution strictly supported on the
boundary of the light cone---only its singular support is the
boundary of the cone. How to compute the action of an operator with
such a singular kernel is much less obvious. The direct summation
along curves provides a fragile, restricted paradigm for curvilinear
integrals, the same way the FFT provides a fragile setting for
shift-invariant problems.

The advantages of our algorithm should now be clear: very general FIOs
can be handled with an asymptotic computational complexity which is
lower than that required for GRT summation, i.e. ($O(N^{2.5} \log N)$
vs. $O(N^3)$), and this without making any curvilinear approximation.
The other argument in favor of the GRT method is the typically low
memory usage. But this equally applies to our method. We will show
that the storage overhead (on top of storing the phase and amplitude)
is negligible and scales like $O(\sqrt{N})$.

We only discussed applications to reflection seismology, but there are
many other areas where nontrivial FIOs are computed routinely, e.g. as
part of solving an inverse problem.  Examples in radar imaging,
ultrasound imaging, and electron microscopy all come to mind. Some
Hough transforms for feature detection in image processing can also be
formulated as FIOs. In short, the ideas presented in this paper may
enable the speed up of fundamental computations in a variety of
problem areas.

\subsection{Related work}

In the case where $\Phi(x,\xi) = x \cdot \xi$, the operator is said to
be pseudodifferential. In this simpler setting, it is known that
separated variables expansions of the symbol $a(x,\xi)$ are good
strategies for reducing complexity. For instance, Bao and Symes
\cite{bao-1996-cpdo} propose a numerical method based on a Fourier
series expansion of the symbol in the angular variable arg $\xi$, and
a polyhomogeneous expansion in $|\xi|$, which is a particularly
effective example of separation of variables.

Another popular approach for compressing operators is to decompose
them in a well-chosen, possibly adaptive basis of $L^2$. Once a sparse
representation is achieved, evaluation simply consists of applying a
sparse matrix in the transformed domain. In the case of 1D oscillatory
integrals, this program was advocated and carried out by Bradie et
al. \cite{bradie-1993-fncoiras} and Averbuch et
al. \cite{averbuch-2000-ecoiamlfb}. In spite of these successes, the
generalization to multiple dimensions has so far remained an open
problem. We will come back to this question in Section \ref{sec:disc},
and in particular discuss the relationship with modern multiscale
transformations such as curvelets \cite{candes-2003-cfio,
  candes-2005-crwpos} and wave atoms \cite{demanet-2006-thesis,
  demanet-2006-wasop}.

We would also like to acknowledge the line of research related to
Filon-type quadratures for oscillatory integrals
\cite{iserles-2004-nqhoift}. When the integrand is of the form $g(x)
e^{ikx}$ with $g$ smooth and $k$ large, it is not always necessary to
sample the integrand at the Nyquist rate. For instance, integration of
a polynomial interpolant of $g$ (Filon quadrature) provides an
accurate approximation to $\int g(x) e^{i k x} \, \d x$ using fewer
and fewer evaluations of the function $g$ as $k \to \infty$. While
these ideas are important, they are not directly applicable in the
case of FIOs. The reasons are threefold. First, we make no notable
assumption on the support of the function to which the operator is
applied, meaning that the oscillations of $\hat{f}(\xi)$ may be on the
same scale as those of the exponential $e^{2\pi i \Phi(x,\xi)}$.
Second the phase does not in general have a simple formula that would
lend itself to precomputations. And third, Filon-type quadratures do
not address the problem of simplifying computations of \emph{several}
such oscillatory integrals at once (i.e. computing a family of
integrals indexed by $x$ in the case of FIOs).

Finally, we remark that FIOs are also interesting when the canonical
relation is nontrivial---that is, multivalued phase---because they
allow to study propagation of singularities of hyperbolic equations in
regimes of multipathing and caustics \cite{hormander-book,
  duistermaat-book}. To mathematicians taking this specialized
viewpoint, the focus of this paper may appear restrictive.  Our
outlook and ambition are different. We find FIOs to be interesting
mathematical objects even when the canonical relation is a graph and
degenerates to the gradient of a phase. Our concern is to understand
their structure from an operational standpoint and exploit it to
design efficient numerical algorithms. In fact, we expect this paper
to be the first of a projected series which will eventually deal with
more complex setups.

\subsection{Contents}

The rest of the paper is organized as follows. Section \ref{sec:analy}
proves all the analytical estimates which support our methodology. In
Section \ref{sec:algo}, we describe algorithms for constructing the
low rank separated approximation, evaluating $(Lf)(x)$, as well as for
evaluating its adjoint, namely, computing $(L^*f)(x)$. Numerical
examples in Section \ref{sec:results} illustrate the properties of our
algorithms. Finally, Section \ref{sec:disc} discusses some related
work and potential alternatives.




\section{Analytical Estimates}
\label{sec:analy}

In this section, we return to a description of the problem in
continuous variables $x$ and $\xi$ to prove estimates on the
separation rank of $e^{2 \pi i \Phi_\ell(x,\xi)}$, where
$\Phi_\ell(x,\xi)$ is the residual phase after linearization about
$\hat \xi_\ell$.

\subsection{Background}

We begin with a lemma which concerns the separation of the exponential
function and whose variations play a central role in modern numerical
analysis.
\begin{lemma}\label{teo:sep_exp}
  Consider the domain defined by $x \in [-A,A]$ for some $A > 0$, and
  $y \in [-1,1]$. For all $\epsilon > 0$ the $\epsilon$-rank
  $r_\epsilon$ of $e^{ixy}$ on $[-A,A] \times [-1,1]$ obeys the
  bound $r_\epsilon \leq r^*_\epsilon$, where
\begin{equation}\label{eq:r_epsilon1}
  r^*_\epsilon = 1 + \max \{ 2 e A \; , \; \log_2(2 \epsilon^{-1}) \}.
\end{equation}
Furthermore, if $A \leq \frac{1}{2 e}$ then the stronger bound
\begin{equation}\label{eq:r_epsilon2}
  r^*_\epsilon = 1 + \frac{\log(2 \epsilon^{-1})}{\log{\frac{1}{e A}}}
\end{equation}
holds as well. In both cases, the corresponding separated
representation is the expansion
\[
|e^{ixy} - \sum_{n = 0}^{r^*_\epsilon - 1} \frac{i^n}{n!} x^n y^n|
\leq \epsilon.
\]
\end{lemma}
\begin{proof}
  The proof is very simple. We start with
\[
\left|e^{i xy} - \sum_{n = 0}^{r - 1} \frac{ (ixy)^n}{n!}\right| \leq
\sum_{n \geq r} \frac{A^n}{n!} \leq \sum_{n \geq r} \left(
  \frac{eA}{n} \right)^n.
\]
It is now straightforward algebra to check that the condition
\[
r \ge \begin{cases} {\log(2 \epsilon^{-1})}/{\log{\frac{1}{e A}}}, & \text{if } eA \leq 1/2, \\
  \max(2 e A, \log_2(2 \epsilon^{-1})), &
  \text{otherwise}, \end{cases}
\]
suffices to bound the right-hand-side by $\epsilon$. Since the
$\epsilon$-rank $r_\epsilon$ is integer-valued, the estimate on $r$
may need to be rounded up to the next integer, hence the precaution of
incrementing the bounds in (\ref{eq:r_epsilon1}) and
(\ref{eq:r_epsilon2}) by one.
\end{proof}

In the next section we will make use of Lemma \ref{teo:sep_exp} to
prove that the nonoscillatory factor $e^{2 \pi i \Phi_\ell(x,\xi)}$
has a separation rank which is independent of $N$. The other factor in
the kernel $a(x,\xi) e^{2\pi i \Phi_\ell(x,\xi)}$, namely, the
amplitude $a(x,\xi)$ is in general a simpler object to study. The
standard assumption in the literature, and also in applications, is to
assume that $a(x,\xi)$ is a smooth symbol of order zero and type
$(1,0)$, meaning that for each pair of integers $(\alpha, \beta)$,
there is a positive constant $C_{\alpha \beta}$ obeying
\[
|\pd^\alpha_\xi \pd^\beta_x a(x,\xi)| \leq C_{\alpha \beta} ( 1 +
|\xi|^2 )^{- |\alpha|/2}. 
\]
For simplicity, we will also assume that $a(x,\xi)$ is compactly
supported in $x$\footnote{This assumption is equivalent to assuming
  that functions in the range of $L$ are themselves compactly
  supported in situations of interest, which ought to be the case for
  accurate numerical computations.}. The nice separation properties of
$a$ are simple consequences of its assumed smoothness.
\begin{lemma}
  Assume $a(x,\xi)$ is a symbol of order zero. Then for all $M > 0$
  there exists $C_M > 0$ such that for all $\epsilon > 0$, the
  $\epsilon$-rank for the separation of $x$ and $\xi$ in $a(x,\xi)$
  obeys
\[
r_\epsilon \leq C_M \, \epsilon^{-1/M}.
\]
\end{lemma}
\begin{proof}
  Perform a Fourier transform of the $C^\infty$, compactly supported
  function $a(\cdot,\xi)$. It suffices to keep $O(\epsilon^{-1/M})$
  Fourier modes to approximate $a(\cdot,\xi)$ to accuracy $\epsilon$
  on its compact support. Each Fourier mode is of the form
  $\hat{a}(\omega,\xi) e^{i \omega x}$, hence separated.
\end{proof}

It goes without saying that the $\epsilon$-rank of the product
$a(x,\xi) e^{2 \pi i \Phi_\ell(x,\xi)}$ is bounded by a constant times
the product of the individual $\epsilon$-ranks, and we now focus on
the real object of interest, the factor $e^{2 \pi i \Phi_\ell}$.

\subsection{Large $N$ asymptotics}\label{sec:estimates_largeN}

In this section we assume that the phase $\Phi(x,\xi)$ is $C^3$ in
$\xi$, only measurable in $x$, and define
\[
C_k = 2 \pi \sup_{x \in [0,1]^2} \sup_{\xi: |\xi| = 1} |\pd_\theta^k
\Phi(x,\xi)| \qquad \mbox{for} \; 0 \leq k \leq 3,
\]
where $\theta = \text{arg} \, \xi$. These constants will enter our
estimates only through the following combinations:
\[
D_2 = C_0 + C_2, \qquad \mbox{and} \qquad D_3 = C_1 + C_3.
\]
As before, we also require homogeneity of order one in $\xi$.
Finally, we let the general angular opening of the cone $W_\ell$ to be
$\frac{2 \alpha}{\sqrt{N}}$ radians, for some constant $\alpha$
(the introduction section proposed $\alpha = \pi$).

The result below is a more precise version of the theorem we
introduced in Section 1.
\begin{theorem}\label{teo:rank1}
 For all $0 < \epsilon \leq 1$, and $N \geq \frac{\alpha^6 D_3^2}{18
    \epsilon^2}$, the
  $\epsilon$-separation rank of $e^{2 \pi i \Phi_\ell(x,\xi)}$ for $x
  \in [0,1]^2$ and $\xi \in W_\ell$ obeys
\begin{equation}\label{eq:rank1_1}
  r_\epsilon \leq 1 + \max \{ \frac{e \sqrt{2}}{2} \alpha^2 D_2 
\; , \; \log_2(4 \epsilon^{-1}) \}.
\end{equation}
Furthermore, if $\alpha$ is admissible in the sense that $\alpha \leq
\sqrt{\frac{\sqrt{2}}{e D_2}}$, then
\begin{equation}\label{eq:rank1_2}
  r_\epsilon \leq 1 + \frac{\log(4 \epsilon^{-1})}{\log{\frac{2 \sqrt{2}}{e \alpha^2 D_2}}}.
\end{equation}
\end{theorem}
\begin{proof} 

  

  Put $r = |\xi|$ and $\theta$ as the angle measured from the vector
  $\xi_\ell$. The phase $\Phi$ can be rewritten as $\Phi(x,\xi) = r
  \phi(x,\theta)$. Let $\xi_1$ be the frequency coordinate along
  $\xi_\ell$ and $\xi_2$ orthogonal to $\xi_1$, so that we can switch
  between polar and Cartesian coordinates using
  \[
  \frac{\pd \Phi}{\pd \xi_1}(x,\xi_\ell) = \phi(x,0), \qquad \mbox{and} \qquad \frac{\pd \Phi}{\pd \xi_2}(x,\xi_\ell) = \phi'(x,0),
  \]
  where the derivative of $\phi$ is taken in $\theta$. The residual
  phase is
  \begin{align*}
    \Phi_\ell(x,\xi) &= \Phi(x,\xi) - \nabla_\xi \Phi(x,\xi_\ell) \cdot \xi \\
    &= r \, \left( \phi(x,\theta) - \cos \theta \phi(x,0) - \sin
      \theta \phi'(x,0) \right).
  \end{align*}
  We can now expand $\phi(x,\theta)$, $\cos \theta$ and $\sin \theta$ in
  a Maclaurin series (around $\theta = 0$) to obtain
  \begin{equation}\label{eq:phir_23}
    \Phi_\ell(x,\xi) = \frac{r \theta^2}{2}
    \left( \phi(x,0) + \phi''(x,0) \right) + 
    \left[ \frac{r \tilde{\theta}^3}{6} \phi'(x,0) + 
      \frac{r \theta^3}{6} \phi'''(x,\overline{\theta}) \right],
  \end{equation}
  for some $\tilde{\theta}$ and $\overline{\theta}$ between $0$ and
  $\theta$ (with $\overline{\theta}$ depending on $x$.)
  
  The $x$ and $\xi$ variables are separated in the first term of
  equation (\ref{eq:phir_23}), so we write
  \[
  f(x) g(\xi) \equiv 2 \pi \left( \phi(x,0) + \phi''(x,0) \right) \frac{r \theta^2}{2}.
  \]
  The term in square brackets is the remainder, and we write
  \[
  R(x,\xi) = 2 \pi \left[ \frac{r \tilde{\theta}^3}{6} \phi'(x,0) +
    \frac{r \theta^3}{6} \phi'''(x,\overline{\theta}) \right].
  \]
  Our strategy will be to choose $N$ large enough so that $R(x,\xi)$
  becomes negligible, hence only the exponential of the first term needs
  to be separated.
  
  Recall that in 2D the frequency domain is the square $[-\frac{N}{2},
  \frac{N}{2}-1]^2$. Since $|\theta| \leq \frac{\alpha}{\sqrt{N}}$ in the
  wedge $W_\ell$, and $r \leq \frac{\sqrt{2}}{2} N$, we have the
  following bounds for the two terms in equation (\ref{eq:phir_23}):
  \[
  |f(x) g(\xi)| \leq \frac{\sqrt{2}}{4} \alpha^2 D_2, \qquad |R(x,\xi)|
  \leq \frac{\sqrt{2}}{12} \frac{\alpha^3}{\sqrt{N}} D_3.
  \]
  It is instructive to notice that the bound on $|fg|$ is independent of
  $N$. That is the reason why we chose the angular opening of the cone
  $W_\ell$ proportional to $N^{-1/2}$ (parabolic scaling).
  
  The first contribution to the separation remainder is given by
  \begin{align*}
    |e^{i (fg + R)} - e^{i fg}| &= |e^{i R} - 1| \\
    &\leq |R| \leq \frac{\sqrt{2}}{12} \frac{\alpha^3}{\sqrt{N}} D_3.
  \end{align*}
  The condition on $N$ ensures precisely that this remainder be
  dominated by $\epsilon/2$.
  
  The second contribution to the total error is due to the separation of
  $e^{i fg}$ itself, and needs to be made smaller than $\epsilon/2$ as
  well. We invoke Lemma \ref{teo:sep_exp} with $f(x) \times
  \sup{|g(\xi)|}$ in place of $x$, $g(\xi)/\sup|g(\xi)|$ in place of
  $y$, and $\epsilon/2$ in place of $\epsilon$. With these choices, $A$
  becomes $\frac{\sqrt{2}}{4} \alpha^2 D_2$, and we obtain the desired
  result.
\end{proof}

\subsection{Small $\epsilon$ asymptotics}

Theorem \ref{teo:rank1} is a special asymptotic result in the case of
large $N$ (problem size) --- or alternatively small $\alpha$ (cone's
angular opening). This regime may not be attained in practice so we
need another result, without restrictions on $N$, and informative for
arbitrarily small $\epsilon$.

To this effect, we need stronger (yet still realistic) smoothness
assumptions on the phase $\Phi$: for each $x$, we require that
$\Phi(x,\xi)$ be a real-analytic function of $\xi$. This condition
implies the bound
\[
2 \pi \sup_{|\xi| = 1} |\pd_\theta^k \Phi(x,\xi)| \leq Q \, k! \, R^{-k},
\]
for some constants $Q$ and $R$. For example, $R$ can be taken as any
number smaller than the uniform radius of convergence in $\theta$, in
which case $Q$ will in general depend on $R$. Let us term such phases,
or functions, \emph{$(Q,R)$-analytic}. As before, we also require
homogeneity in $\xi$.

\begin{theorem}\label{teo:rank2}
  Assume $\Phi_\ell(x,\xi)$ is measurable in $x$, and $(Q,R)$-analytic
  in $\xi$, for some constants $Q$ and $R$. Assume that $\alpha$ is
  admissible in the sense that
  \[
  \alpha < \min \{ \, \frac{r \sqrt{N}}{2}, \frac{R}{\sqrt{\sqrt{2} Q}} \}.
  \]
  Then for all $0 < \epsilon \leq 1$, the $\epsilon$-separation rank
  of $e^{2 \pi i \Phi_\ell(x,\xi)}$ for $x \in [0,1]^2$ and $\xi \in
  W_\ell$ obeys
  \[
  r_\eps \le C_p \, \epsilon^{-p}, \quad \forall p : p > \frac{2}{\log_2
    \left( \frac{R \sqrt{N}}{\alpha} \right)}.
  \]
\end{theorem}

\begin{proof}
  Throughout the proof, $x \in [0,1]^2$ and $\xi \in W_\ell$.  Using
  the smoothness assumption on $\Phi_\ell$, we can repeat the
  reasoning of the proof of Theorem \ref{teo:rank1} and obtain the
  convergent series
  \[
  2 \pi \Phi_\ell(x,\xi) = \sum_{k=0}^\infty f_k(x) g_k(\xi),
  \]
  where $f_k(x) = 2 \pi \phi^{(k)}(x,0)$ (the differentiations are in
  $\theta$) and
\[
g_0(\xi) = r (1 - \cos \theta), \qquad g_1(\xi) = r (\theta - \sin
\theta), \qquad g_k(\xi) = \frac{r \theta^k}{k!}.
\]
We denote the bound $|f_k(x) g_k(\xi)| \leq A_k$, with
\[
A_0 = \frac{\sqrt{2}}{4} Q \alpha^2, \qquad A_1 = \frac{\sqrt{2}}{12}
Q \frac{\alpha^3}{R \sqrt{N}}, \qquad A_k = \frac{\sqrt{2}}{2} Q N
\left( \frac{\alpha}{R \sqrt{N}} \right)^k \quad \mbox{for } k > 2.
\]

Our strategy will be to call upon Lemma \ref{teo:sep_exp} for the
first few factors $e^{i f_k g_k}$, in order to obtain a separation
rank $r_k$ and an error $\epsilon_k$ for each of them:
\begin{equation}\label{eq:sep_fkgk}
  \left| \, e^{i f_k g_k} - \sum_{n=0}^{r_k - 1} \frac{i^n}{n !} f_k^n(x) g^n_k(\xi) \right| \leq \epsilon_k.
\end{equation}
We will perform this operation for each $k < K$, with $K$ large
enough, to be determined. Once the separation of each factor is
available, we can write
\[
e^{i \sum_{k = 0}^{K-1} f_k g_k} = \prod_{k=0}^{K-1} e^{i f_k g_k},
\]
and obtain the bound on the overall separation rank as the
\emph{product} $\prod_{k=0}^{K-1} r_k$.

There are two sources of errors we must contend with:
\begin{itemize}
\item \emph{Truncation in $k$}. The factors $e^{i f_k g_k}$ for $k
  \geq K$, will be deemed negligible if their combined contribution
  results in an overall error smaller than $\epsilon/2$, meaning
\begin{equation}\label{eq:truncink}
| e^{2 \pi i \Phi_\ell} - e^{i \sum_{k = 0}^{K-1} f_k g_k} | \leq \frac{\epsilon}{2}.
\end{equation}
The left hand side is bounded by $|\sum_{K}^\infty f_k g_k|$. Using
the bound we stated earlier on $A_k$, and the admissibility condition on $\alpha$, a bit of algebra shows that (\ref{eq:truncink}) is satisfied for
\begin{equation}\label{eq:K}
  K = \lceil \frac{\log \left( 2 \sqrt{2} QN \epsilon^{-1} \right)} {\log \left( \frac{R
        \sqrt{N}}{\alpha} \right)} \rceil
\end{equation}
(meaning the smallest integer greater than the quotient inside the brackets). This quantity in turn obeys $\log_2(8 \epsilon^{-1}) \leq K \leq \log_2(16 \epsilon^{-1}) $.

\item \emph{Truncation in $n$}. The truncation errors from
  (\ref{eq:sep_fkgk}) must be made sufficiently small so that their
  combined contribution also results in an overall error smaller than
  $\epsilon/2$, meaning
\begin{equation}\label{eq:truncn}
  \left| \, \prod_{k=0}^{K-1} e^{i f_k g_k} - \prod_{k=0}^{K-1} \sum_{n=0}^{r_k -
      1} \frac{i^n}{n !} f_k^n(x) g^n_k(\xi) \, \right| \leq \frac{\epsilon}{2}.
\end{equation}

Easy manipulations\footnote{To justify this step, put $E_k(x,\xi) =
  e^{i f_k(x) g_k(\xi)}$ and start from the identity
\begin{align*}
  \prod_{k = 0}^{K-1} (E_k + \epsilon_k) - \prod_{k = 0}^{K-1} E_k &=
  \sum_j \epsilon_j \prod_{k \ne
    j} (E_k + \tau_{jk} \epsilon_k) \\
  &= \sum_j \epsilon_j E^{-1}_j \prod_{k = 0}^{K-1} (E_k + \tau_{jk}
  \epsilon_k)
\end{align*}
where $\tau_{jk} = 0$ if $j \leq k$, and $\tau_{jk} = 1$ if $j >
k$. Then make use of the bound $(1 + \frac{\epsilon}{3K})^K <
e^{\epsilon/3} \leq e^{1/3} < 3/2$.} show that (\ref{eq:truncn})
follows from the bound
\[
\epsilon_k = \frac{\epsilon}{3 K}.
 \]
 (Recall that $K$ is comparable to $\log (C \epsilon^{-1})$.)

 Such a bound holds if, in turn, we take $r_k$ large enough. The
 admissibility condition on $\alpha$ ensures, among others, that we
 can invoke the strong version of Lemma \ref{teo:sep_exp}, namely
 equation (\ref{eq:r_epsilon2}), and obtain
\begin{equation}\label{eq:rk}
r_k \leq 1 + \frac{\log(2 \epsilon_k^{-1})}{\log{\frac{1}{e A_k}}}.
\end{equation}

\end{itemize}

It now remains to estimate $\prod_{k=0}^{K-1} r_k$, where $r_k$ is
given by equation (\ref{eq:rk}) and $K$ by equation (\ref{eq:K}). We
treat the first two factors independently: we can check from the
bounds on $A_0$ and $A_1$, and the admissibility condition on
$\alpha$, that
\[
r_0 r_1 \leq C \log \left[ 2 \epsilon^{-1} \log (2 \epsilon^{-1}) \right].
\]
As for the case $k \geq 2$,
\begin{align*}
  r_k &\leq 1 + \frac{\log(6 K \epsilon^{-1})} {\log\left[
      \frac{\sqrt{2}}{Q N}
      \left( \frac{R \sqrt{N}}{\alpha} \right)^{k} \right]} \\
  &\leq \frac{\log(C \epsilon^{-1} \log (2 \epsilon^{-1})) + k \log
    \left( \frac{R \sqrt{N}}{\alpha} \right) }{\log(C) +
    k \log \left( \frac{R \sqrt{N}}{\alpha} \right)} \\
  &\equiv \frac{A + k}{B+k}. \qquad \qquad (k \geq 2)
\end{align*}
We only simplified notations in the last line. Notice that $A > B$,
and that $B + k \geq 1$ when $k \geq 2$. We will assume without loss
of generality that $A$ and $B$ are integers.The value of the product
$\prod r_k$ can only increase if we replace the initial bound $0\leq k
< K$, by the condition that the bound on $r_k$ be greater than $2$. So
we certainly have
\begin{align*}
  r_\epsilon &\leq \prod_{k \geq 2 : r_k \geq 2} r_k 
\leq \frac{A+2}{B+2} \; \frac{A+3}{B+3} \; \ldots 
\; \frac{A+A}{B+A} \\
  &= \frac{(2A)! / (A+1)!}{(B+A)! / (B+1)!}. 
\end{align*}

We can now make use of the two-sided Stirling bound
\[
\sqrt{2 \pi} \, n^{n+1/2} e^{-n + \frac{12}{n+1}} \leq n! \leq \sqrt{2
  \pi} \, n^{n+1/2} e^{-n + \frac{12}{n}}
\]
to obtain
\begin{align*}
  r_\epsilon &\leq C \, \frac{(2A)^{2A} (A+1)^{-(A+1)}}{(A+B)^{A+B} (B+1)^{-(B+1)}} \\
  &\leq C \, 2^{2A} \, \frac{A^{2A}}{(A+1)^{(A+1)} (A+B)^{A-1}} \, \frac{(B+1)^{(B+1)}}{(A+B)^{B+1}} \\
  &\leq C \, 2^{2A}.
\end{align*}
In turn,
\[
2^{2A} \leq \left( C \epsilon^{-1} \log (2 \epsilon^{-1}) \right)^{\frac{2}{\log_2 \left( \frac{R \sqrt{N}}{\alpha} \right) }},
\]
which concludes the proof.
\end{proof}

The lower the fractional exponent of $\epsilon^{-1}$ the faster the
convergence of separated expansions. Theorem \ref{teo:rank2} shows
exactly which factors can make this exponent arbitrarily small:
\begin{itemize}
\item large grid size $N$, or
\item small angular opening constant $\alpha$, or
\item large radius of analyticity $R$ of the phase in arg $\xi$
  (uniformly in $x$).
\end{itemize}
Observe that the rank bound \emph{decreases} as $N$ increases.

Theorem \ref{teo:rank2} assumes that the residual phase function
$\Phi_\ell(x,\xi)$ is $(Q,R)$-analytic in $\xi$. The variation below
follows the same path of reasoning, and is useful when
$\Phi_\ell(x,\xi)$ is only $C^\infty$ in $\xi$ for $\xi\not=0$.


\begin{theorem}\label{teo:lexing}
  Assume $\Phi_\ell(x,\xi)$ is $C^\infty$ in $\xi$ for $\xi\not=0$.
  For any $p > 0$, there exists two constants $C_p$ and $C_p'$ such
  that for any $N$, the $\epsilon$-separation rank with $\eps=C_p \, N^{-p}$
  is bounded by $C_p' \, \log N$.
\end{theorem}
\begin{proof}
  The structure of the proof is similar to that of Theorem
  \ref{teo:rank2}. One only needs to keep the first $2p+2$ term of the
  series
  \[
  2 \pi \Phi_\ell(x,\xi) = \sum_{k=0}^\infty f_k(x) g_k(\xi).
  \]
  in order to have $\eps=C_p \, N^{-p}$ for some constant $C_p$ which depends
  only on $p$ and $\Phi_\ell$.
  The product $\prod_{k=0}^{2p+1} r_k$ upper bounds the overall
  separation rank, and is less than $C_p' \, \log N$ for some
  constant $C_p'$ which only depends on $p$.
\end{proof}

In many computational problems, the mesh size $N^{-1}$ is linked
directly to the desired accuracy $\eps$, usually in the form of a
power law, e.g.~$\eps = O(N^{-p})$ for some constant $p$. Therefore,
Theorem \ref{teo:lexing} is interesting for practical reasons.




\section{Algorithm}
\label{sec:algo}

For notational convenience, we assume in this section that the
amplitude is identically equal to one; that is, we focus on the
so-called (discretized) Egorov operator
\begin{equation}\label{eq:Egorov}
  (L f)(x) = \frac{1}{N} \sum_{\xi  \in \Omega} e^{2 \pi i \Phi(x,\xi)} \hat{f}(\xi).
\end{equation}
Both in practice (Section \ref{sec:results}) and in theory (Section
\ref{sec:analy}), one can easily take care of general amplitude terms.

The algorithm for computing \eqref{eq:Egorov} has two main components:
\begin{itemize}
\item {\em Preprocessing step}. Given the residual phase $\Phi_\ell(x,\xi)
  \equiv \Phi(x,\xi)- \Grad_\xi \Phi(x,\hat{\xi}_\ell)\cdot\xi$, this
  step constructs, for each wedge $W_\ell$, a low rank separated
  approximation
  \[
  \left| \,
  e^{2\pi i \Phi_\ell(x,\xi)} -
  \sum_{t=1}^{q} \gamma^x_{\ell t}(x) \gamma^\xi_{\ell t}(\xi)
  \, \right|
  \le \eps.
  \]
  The functions $\{\gamma^x_{\ell t}(x)\}$ and $\{\gamma^\xi_{\ell
    t}(\xi)\}$, or their compressed versions, are then stored for use
  in the next step.
\item {\em Evaluation step}. Given a function $f$, this step computes
  $(Lf)(x)$ approximately by 
  \[
  (Lf)(x) \approx \frac{1}{N} \sum_\ell \sum_t \gamma^x_{\ell t}(x)
  \sum_\xi e^{2\pi i \Grad_\xi \Phi(x,\hxw)\cdot\xi} \left[ \,
    \gamma^\xi_{\ell t}(\xi) \chi_\ell (\xi) \hat{f}(\xi) \, \right].
  \]
\end{itemize}
The preprocessing step is performed only once for a fixed phase
function $\Phi(x,\xi)$. The family of functions $\{ \gamma^x_{\ell
  t}(x) \}$ and $\{ \gamma^\xi_{\ell t}(\xi) \}$ should of course be
used again and again to compute $(Lf)(x)$ for different inputs $f$.

In Sections \ref{sec:algo-dete} and \ref{sec:algo-rand}, we propose
two different approaches for constructing the families $\{
\gamma^x_{\ell t}(x) \}$ and $\{ \gamma^\xi_{\ell t}(\xi) \}$. Section
\ref{sec:algo-eval} describes the details of the evaluation step.
Finally, Section \ref{sec:algo-evaladj} outlines the algorithm for
rapidly applying the adjoint operator. In this section, we calculate
time and storage complexity under the assumption of large grids,
i.e. that of Theorem \ref{teo:rank2}. For other kinds of asymptotics,
one may need to adjust these estimates with a multiplicative $\log N$
factor, which is typically negligible.
\subsection{Preprocessing step: deterministic approach}
\label{sec:algo-dete}

We first describe a deterministic approach for constructing the low
rank separated expansion, based on a Taylor expansion, exactly as in
the proof of Lemma \ref{teo:sep_exp}. For each wedge $W_\ell$, the
strategy consists of the following sequence of steps:
\begin{enumerate}
\item Construct a low rank separated approximation of
  $\Phi_\ell(x,\xi)$. This is done by truncating the polar coordinates
  Taylor expansion to the $(2p+1)$st term
  \[
  \Phi_\ell(x,\xi) \approx |\xi| \, 
  \sum_{k=1}^{2p+1} c_{\ell k}(x) (\theta-\tw)^k. 
  \]
  Here $p$ is a constant that determines the level of accuracy.
\item For each $k$ construct a separated expansion of $e^{2\pi i
    c_{\ell k}(x) \, |\xi| (\theta-\tw)^k}$. This is done by
  truncating the Taylor expansion to the first $d_{\ell k}$ terms
  \[
  e^{2\pi i c_{\ell k}(x) \, |\xi|(\theta-\tw)^k} \approx
  \sum_{m=0}^{d_{\ell k}-1} \beta^x_{\ell k m}(x) \beta^\xi_{\ell k
    m}(\xi).
  \]
  The value of each $d_{\ell k}$ is also chosen to obtain a good accuracy.
\item Combine the separated expansions for $k =1, \ldots, 2p+1$ into
  one separated representation for $e^{2\pi i \Phi_\ell(x,\xi)}$.
  Simply expanding the product of the expansions obtained in the
  previous step would be sufficient for proving a theorem like those
  presented in Section \ref{sec:analy} but in practice though, the
  number of terms in the expansion is too large and far from optimal.
  We thus combine the product of separated expansions two-by-two with
  the compression procedure to be described next, and repeat the
  process until there is only one separated expansion left. The final
  expansion provides us with the required functions $\{ \gamma^x_{\ell
    t}(x) \}$ and $\{ \gamma^\xi_{\ell t}(\xi) \}$.
\end{enumerate}

The compression procedure used to combine the product of two separated
expansions is quite standard. Suppose we only have two expansions (the
subscript $\ell$ is implicit) and write their product as
\begin{multline*}
  \left(\sum_{m_1=0}^{d_{1}-1} \beta^x_{1m_1}(x)
    \beta^{\xi}_{1m_1}(\xi) \right)
  \left( \sum_{m_2=0}^{d_{2}-1} \beta^x_{2m_2}(x) \beta^{\xi}_{2m_2}(\xi) \right) \\
  = \sum_{m_1,m_2} \left( \beta^x_{1m_1}(x) \beta^x_{2m_2}(x)
  \right)
  \left( \beta^{\xi}_{1m_1}(\xi)  \beta^{\xi}_{2m_2}(\xi) \right)
   := \sum_{m} c^x_{m}(x) c^{\xi}_{m}(\xi).
\end{multline*}
We adopt the matrix notation and introduce
\[
(A)_{x,m} = c^x_m(x), \quad (B^*)_{m,\xi} = c^\xi_m(\xi).
\]
The problem is to find two matrices $\tilde{A}$ and $\tilde{B}$ which
have far fewer columns than $A$ and $B$, and yet obeying $\tilde{A}
\tilde{B}^* \approx AB^*$. This may be achieved by means of the QR
factorization and of the SVD:
\begin{enumerate}
\item Construct QR factorizations $A= Q_A R_A$ and $B = Q_B R_B$. 
\item Compute the singular value decomposition of $R_A R_B^*$ and
  truncate the singular values below a threshold $\eps$ together with
  their associated left and right singular vectors, i.e. $R_A R_B^*
  \approx U_M S_M V_M^*$ where $S_M$ is a truncated diagonal matrix of
  singular values.
\item Set $\tilde{A} = Q_A U_M S_M$ and $\tilde{B} = Q_B V_M$.
\end{enumerate}

Suppose $A$ is $m \times q$ and $B$ is $n \times q$ with both $m$ and
$n$ much larger than $q$. The computational complexity of the
compression procedure is $O((m+n) q^2)$. In our setup, $m = |X| =
N^2$, $n = |W_\ell| = O(N^{1.5})$, and $q$, the rank bound, is
uniformly bounded in $N$ (Theorem \ref{teo:rank2} shows that $q$ is
bounded by a small fractional power of $\eps$, independently of $N$).
Therefore, the complexity of a single compression procedure is
$O(N^2)$. Since this needs to be carried out $2p-1$ times for each of
the $\sqrt{N}$ wedges, the overall complexity of the deterministic
preprocessing is $O(\sqrt{N} \times N^2) = O(N^{2.5})$ where the
constant is directly related to the rank bounds of Section
\ref{sec:analy}.

Next, let us consider the storage requirement. For each wedge, the
size of the final separated expansion is $O(N^2)$. Since there are
$\sqrt{N}$ wedges, the total storage requirement is $O(N^{2.5})$,
which can be costly when $N$ is large. For example, in a typical
problem with $N = 1024$ and $q =20$, the total storage would be about
10 GB assuming double precision is used. Our second approach to solve
the preprocessing step addresses this issue and requires dramatically
less storage space.

\subsection{Preprocessing step: randomized approach}
\label{sec:algo-rand}

This section describes a randomized approach for computing the
functions $\{ \gamma^x_{\ell t}(x) \}$ and $\{ \gamma^\xi_{\ell
  t}(\xi) \}$ for a fixed $\ell$. The method is based on the work
presented in Kapur and Long \cite{kapur-1997-IES3}. We use matrix
notations and set $A$ to be the matrix defined by
\begin{equation}
  \label{eq:probform}
  A_{x,\xi} := e^{2\pi i \Phi_\ell(x,\xi)}, \quad x\in X, \xi \in W_\ell. 
\end{equation}
The matrix $A$ is $m$ by $n$ with $m = N^2$ and $n = O(N^{1.5})$.
Assume the prescribed error $\eps$ is fixed, Theorem \ref{teo:rank2}
tells us that there exists a low rank factorization of $A$ with rank
$r_\eps = O(1)$ (again, by this we mean that $r_\eps$ is bounded by a
constant independent of $N$, although not independent of
$\epsilon$). Using this knowledge, the following randomized method
finds an approximate factorization
\[
A \approx U T, 
\]
where $U$ is of size $m \times q $, $T$ is $q \times n$ and $q=O(1)$
in $N$.
\begin{enumerate}
\item Select a set $C$ of $r$ columns taken from $A$ uniformly at
  random, and define $A_{[C]}$ to be the
  submatrix formed by these columns. In practice, a safe choice is to
  take $r$ about three times larger than the (unknown) $r_\eps$.

\item Compute the singular value decomposition $A_{[C]} \approx U S V^*$
  where the diagonal of $S$ contains only the singular values greater
  than the threshold $\eps$. Since $A$ has a separation rank $r_\eps =
  O(1)$, we expect $U$ to be of size $m\times q$ where $q$ is about
  $r_\epsilon$.

\item Select a set $R$ of $r$ rows taken from $A$ uniformly at random,
  and define $A_{[R]}$ to be the submatrix formed by these rows. Similarly, 
  let $U_{[R]}$ be the submatrix of $U$ containing the same rows.

\item Set $T = U^+_{[R]} A_{[R]}$ where $U^+_{[R]}$ is the pseudo-inverse
  of $U_{[R]}$.

\item The matrices $U$ and $T$ provide an approximate factorization,
  i.e. $A \approx U T$. We identify the columns of $U$ with the family
  $\{ \gamma^x_{\ell t}(x) \}$, and the rows of $T$ with
  $\{\gamma^\xi_{\ell t}(\xi) \}$.
\end{enumerate}

This randomized approach works well in practice although we are not
able to offer a rigorous proof of its accuracy, and expect one to be
non-trivial. We merely argue that the validity of this methodology
hinges on the following observations:
\begin{itemize}
\item First, the columns of $A$ are highly correlated. Following the
  arguments in Section \ref{sec:analy}, it is not difficult to show
  that a pair of columns with nearby values of the frequency index
  $\xi \in W_\ell$ have a large inner product. Therefore, as we sample
  uniformly at random, we get a good coverage of the set $W_\ell$
  (leaving no large hole) and as a result, the sampled columns nearly
  span the space generated by the columns of $A$.  Note that one could
  also use a deterministic regular sampling strategy; for instance, we
  could take a Cartesian subgrid as a subset of $W_\ell$.  We observed
  that in practice, the probabilistic approach provides slightly
  better approximations.
  
  

  
\item As the SVD routine is numerically stable, it allows us to extract
  an orthobasis of the column space of $A_{[C]}$ in a robust way.
\item By construction, the columns of $U$ are orthonormal. Results
  from random projection and the geometry of high-dimensional spaces
  imply that, as long as $U$ does not correlate with the canonical
  orthobasis, the columns of $U_{[R]}$ are almost orthogonal as well. This
  allows us to recover the matrix $T$ in a stable and robust fashion.
\end{itemize}

The computational complexity of this randomized approach is quite low.
The SVD step has a complexity of $O(m r^2) = O(N^2)$, while the matrix
product $T = U^+_{[R]} A_{[R]}$ takes $O(n r q) = O(N^{1.5})$
operations. Therefore, for each $\ell$, the complexity of the
randomized approach is $O(N^2)$. Since the same procedure needs to be
carried out for all the $\sqrt{N}$ wedges, the overall complexity is
$O(N^{2.5})$.

Often we do not know the exact value of $r_\eps$. Instead of setting
$r$ conservatively to be an unnecessarily large number, this
difficulty is addressed as follows: we begin with a small $r$, and
check whether $q$ is significantly smaller than $r$. If this is the
case, we accept the factorization. Otherwise, we double $r$ and
restart the process. Geometrical increase guarantees that the work
wasted (due to unsuccessful attempts) is bounded by the work of the
final successful attempt. In practice, we accept the result when $q
\le r/3$, and this criterion seems to work well in our numerical
experiments. A more conservative test certainly improves the
reliability of the factorization but increases the running time.
 
We finally examine the storage requirement. A naive approach is to
store the matrices $U$ and $T$ for each wedge $W_\ell$. As $T$ is much
smaller than $U$ in size, the storage requirement for each wedge is
roughly the size of $U$, which is $N^2 q = O(N^2)$. Multiplying this
by the number of wedges gives a total storage requirement of
$O(N^{2.5})$, which can be quite costly for large $N$ as already
mentioned in the last section. We propose to store the matrices $V
S^{-1}$ and $U^+_{[R]}$ instead. Both matrices only require storage of
size $O(rq) = O(1)$. Whenever we need $U$ and $T$, we form the
products $U = A_{[C]} V S^{-1}$ and $T = U^+_{[R]} A_{[R]}$. Note that the
elements of the matrices $A_{[C]}$ and $A_{[R]}$ are given explicitly by the
formula \eqref{eq:probform} and there is of course no need to store
them at all. Putting it differently, we rewrite the computed
factorization as
\begin{equation}
  A \approx A_{[C]} \, V S^{-1}  U^+_{[R]} \, A_{[R]} 
  \label{eq:probnewfac}
\end{equation}
and store only the matrices $V S^{-1}$ and $U^+_{[R]}$.

We would like to point out that such a scheme is not likely to work
for the deterministic approach. The main reason is that the
deterministic approach involves multiple compression procedures which
make use of QR factorizations and SVD decompositions. These numerical
linear algebra routines are quite complicated, and therefore, it would
be difficult to relate the resulting low-rank factorization with the
elements of the matrix $A$, which have the simple form
\eqref{eq:probform}.

\subsection{Comparison}

Table \ref{tbl:cmp} compares the deterministic and the randomized
approaches in view of the computational complexity and storage
requirement. The deterministic approach has the advantage of
guaranteeing an accurate low rank separation. However, the constant in
the time complexity can be quite large as for each wedge, it requires
$2p$ compression procedures to combine multiple separated expansions
into a single one. Moreover, since the compression step uses QR
factorizations and SVDs, we are forced to store the final expansion,
which can be quite costly for large $N$. In practice, the randomized
approach constructs a near optimal low rank expansion with very high
probability, requires very low storage space, and enjoys a
significantly lower constant in time complexity since it does not
utilize repeated QR factorizations or singular value decompositions.

\begin{table}[htb]
  \begin{center}
        \begin{tabular}{|c|cc|}
          \hline
          & time & storage\\
          \hline
          randomized & $O(N^{2.5})$ (small constant) & $O(\sqrt{N}) $\\
          deterministic& $O(N^{2.5})$ (large constant) & $O(N^{2.5}) $\\ 
          \hline
        \end{tabular}
  \end{center}  
  \caption{Comparison of the deterministic and randomized approaches.}
  \label{tbl:cmp}
\end{table}

\subsection{Evaluation step}
\label{sec:algo-eval}

Once the families $\{ \gamma^x_{\ell t}(x) \}$ and $\{
\gamma^\xi_{\ell t}(\xi) \}$ are available, we use the approximation
\[
(Lf)(x) \approx \frac{1}{N} \sum_\ell \sum_t \gamma^x_{\ell t}(x)
\sum_\xi e^{2\pi i \Grad_\xi \Phi(x,\hxw)\cdot\xi} \left[ \,
  \gamma^\xi_{\ell t}(\xi) \chi_\ell(\xi) \hat{f}(\xi) \, \right]
\]
to evaluate $L f(x)$. The algorithm simply carries out the evaluation
step by step:
\begin{enumerate}
\item Compute $\hat{f}$, the Fourier transform of $f$. 
\item For each $\ell$ and $t$, form $\hat{f}_{\ell t}(\xi) :=
  \gamma^\xi_{\ell t}(\xi) \chi_\ell(\xi) \hat{f}(\xi)$. 
\item For each $\ell$ and $t$, compute $g_{\ell t}(x) := \sum_\xi
  e^{2\pi i \Grad_\xi \Phi(x,\hxw)\cdot\xi} \hat{f}_{\ell t}(\xi)$. 
\item Compute $(Lf)(x) \approx \frac{1}{N} \sum_\ell \sum_t
  \gamma^x_{\ell t}(x) g_{\ell t}(x)$. 
\end{enumerate}

The only step that requires attention is the third: it asks to
evaluate the Fourier series $\sum_\xi e^{2\pi i \Grad_\xi
  \Phi(x,\hxw)\cdot\xi} \hat{f}_{\ell t}(\xi)$ at the $N^2$ points $\{
\Grad_\xi \Phi(x,\hxw): x\in X\}$. Even though $X$ is a Cartesian
grid, the warped grid $\{ \nabla_\xi \Phi(x,\hat{\xi}_\ell): x\in X\}$
is no longer so. In fact, the formula for $g_{\ell t}$ is a nonuniform
Fourier transform of the second kind, a subject of considerable
attention
\cite{anderson-1996-rcdft,beylkin-1995-fftfws,greengard-2004-anfft,nguyen-1999-rfmnfft,potts-2001-fftndt}
since the seminal paper of Dutt and Rokhlin \cite{dutt-1993-fftnd}.
We adopt the approach introduced in the latter paper, and following
their notations, set
\begin{itemize}
\item
  $m=4$, $q=8$ and $b=0.425$ for 6 digits of accuracy,
\item
  $m=4$, $q=16$ and $b=0.785$ for 11 digits of accuracy.
\end{itemize}
We specify these parameter values because they impact the numerical
accuracies we will report in the next section, and because it will
help anyone interested in reproducing our results.

The algorithm in \cite{dutt-1993-fftnd} generally assumes that the
Fourier coefficients are supported on the full grid $\Omega$ which is
symmetric with respect to the origin. For each $\ell$, the support of
$\hat{f}_{\ell t}(\xi)$ is $W_\ell$, which is to say that most of the
values of the input on the grid $\Omega$ are zero.  To speed up the
nonuniform fast Fourier transform, each wedge $W_\ell$, which is close
to either one of the diagonals, is sheared by 45 degrees so that it
becomes approximately horizontal or vertical. Notice that 45 degree
shearing of $\hat{f}_{\ell t}(\xi)$ is a simple relabeling of the
array. In addition, all wedges are then translated so that their
support fits in a rectangle of smaller volume centered around the
origin. As the nonuniform FFT \cite{dutt-1993-fftnd} asks to compute
the FFT of the input data (and then finds a way of interpolating the
result on an untructured grid), we gain efficiency since the input
array is now of smaller size. Mathematically, the shearing operation takes
the form
\[
\xi' = M \xi - \xi_c, 
\]
where $M$ is either the identity or a 45-degree shear matrix and
$\xi_c$ is a translation parameter. Thus, we organize the computations
as in 
\begin{align*}
  \sum_\xi e^{2\pi i \Grad_\xi \Phi(x,\hxw)\cdot\xi} \hat{f}_{\ell
    t}(\xi)
  & = \sum_{\xi'} e^{2\pi i \Grad_\xi \Phi(x,\hxw)\cdot M^{-1} (\xi' + \xi_c)}  \hat{f}_{\ell t}(M^{-1} (\xi'+\xi_c))\\
  & = e^{2\pi i \Grad_\xi \Phi(x,\hxw)\cdot M^{-1} \xi_c} \,
  \sum_{\xi'} e^{2\pi i \Grad_\xi \Phi(x,\hxw)\cdot M^{-1} \xi'}
  \hat{f}_{\ell t}(M^{-1} (\xi' + \xi_c)),
\end{align*}
where the final summation is a nonuniform Fourier transform at points
$(M^*)^{-1} \Grad_\xi \Phi(x,\hxw)$. In condensed form, the
oscillatory modes of the function we wish to evaluate are centered
around a center frequency; we factor out this frequency, interpolate
the residual, and add the factor back in; for the same accuracy,
interpolating the smoother residual requires a smaller computational
effort.

A two-dimensional nonuniform fast Fourier transform takes $O(N^2 \log
N)$ operations. This operation needs to be repeated $q = O(1)$ times
for each one of the $\sqrt{N}$ wedges. Therefore, the overall
complexity is $O(N^{2.5} \log N)$.


\subsection{Evaluating the adjoint operator}
\label{sec:algo-evaladj}

We conclude this section by presenting how to rapidly apply the
adjoint Fourier integral operator.  Begin by expanding the Fourier
transform in \eqref{eq:FIO} and write
\[
(Lf)(x) = \int \left( \int e^{2\pi i (\Phi(x,\xi)-y\cdot\xi)}
  \d\xi\right) f(y) \d y.
\]
for $x,y,\xi \in \R^2$. The adjoint operator is then given by
\begin{eqnarray*}
  (L^*f)(x) &=& \int 
\left( \int e^{- 2\pi i (\Phi(y,\xi)-x\cdot\xi)} \d\xi\right) f(y) \d y\\
  &=& \int 
\left( \int e^{- 2\pi i \Phi(y,\xi)} f(y) \d y \right) e^{ 2\pi i x\cdot\xi} \d\xi 
\end{eqnarray*}
or equivalently as 
\[
\widehat{ (L^*f) }(\xi) = \int e^{- 2\pi i \Phi(y,\xi)} f(y) \d y
\]
in the Fourier domain.  Similarly, one readily checks that the adjoint
of the discrete-time FIO is given by the formula
\[
\widehat{ (L^*f) }(\xi) = \frac{1}{N} \sum_y e^{- 2\pi i \Phi(y,\xi)} f(y), 
\]
where $\xi\in \Omega$ and $y \in X$.

Now follow the same set of ideas as in Section \ref{sec:algo-eval}, and
decompose $L^*$ as
\begin{eqnarray*}
  \widehat{ (L^*f) }(\xi) &=& \frac{1}{N} \sum_\ell \chi_\ell(\xi) \sum_y e^{- 2\pi i \Phi(y,\xi)} f(y)\\
  &=&       \frac{1}{N} \sum_\ell \chi_\ell(\xi) \sum_y e^{-2\pi i \Phi_\xi(y,\hxw)\cdot\xi} e^{-2\pi i \Phi_\ell(y,\xi)} f(y)\\
  &=& \frac{1}{N} \sum_\ell \chi_\ell(\xi) \sum_y e^{-2\pi i \Phi_\xi(y,\hxw)\cdot\xi} \sum_t \overline{\gamma^x_{\ell t}(y) \gamma^\xi_{\ell t}(\xi)}  f(y)\\
  &=&       \frac{1}{N} \sum_\ell \sum_t \chi_\ell(\xi) \overline{\gamma^\xi_{\ell t}(\xi)} \sum_y e^{-2\pi i \Phi_\xi(y,\hxw)\cdot\xi} \left( \overline{\gamma^x_{\ell t}(y)} f(y) \right). 
\end{eqnarray*}
The right-hand side of the last equation provides the key steps of the
algorithm. 
\begin{enumerate}
\item 
  For each $\ell$ and $t \leq q$, compute $f_{\ell t}(y) := \overline{\gamma^x_{\ell t}(y)} f(y)$. 
\item For each $\ell$ and $t \leq q$, compute $g_{\ell t}(\xi) :=
  \sum_y e^{-2\pi i \Phi_\xi(y,\hxw)\cdot\xi} f_{\ell t}(y)$ using the
  nonuniform fast Fourier transform of the first kind, see
  \cite{dutt-1993-fftnd,greengard-2004-anfft} for details.
\item Compute $\widehat{ (L^*f) }(\xi) \approx \frac{1}{N} \sum_\ell
  \sum_t \chi_\ell(\xi) \overline{\gamma^\xi_{\ell t}(\xi)} g_{\ell
    t}(\xi)$.
\item Finally, take an inverse 2D FFT to get $(L^*f)(x)$.
\end{enumerate}
Clearly, all the results and discussions concerning the matrix vector
product $L f$ apply here as well.

\section{Numerical Results}
\label{sec:results}

This section presents several numerical examples to demonstrate the
effectiveness of the algorithms introduced in Section
\ref{sec:algo}. Our implementation is in Matlab and all the
computational results we are about to report were obtained on a
desktop computer with a 2.6 GHz CPU and 3 GB of memory. We have
implemented both the deterministic and randomized approaches for the
preprocessing step. We choose to report the timing and accuracy
results of the randomized approach only since it requires less time
and storage as shown in Section \ref{sec:algo-rand}.

We first study the error of the separated approximation generated by
the randomized preprocessing step. For $x=(x_1,x_2)$ and
$\xi=(\xi_1,\xi_2)$, set the phase function to be
\begin{equation}
\Phi_{\pm}(x,\xi) =  x\cdot\xi \pm \sqrt{ 
  r^2_1(x)\xi_1^2 + r^2_2(x)\xi_2^2}.
\label{eq:nrphase}
\end{equation}
We show in the Appendix that the transformation, which for each $x$
integrates $f$ along an ellipse centered at $x$ and with axes of
length $r_1(x)$ and $r_2(x)$, can be cast as a sum $L_+ + L_{-}$ of
two FIOs given by
\begin{equation}\label{eq:ellipses_FIO}
(L_\pm f)(x) = \int a_{\pm}(x,\xi) e^{2\pi i \Phi_{\pm}(x,\xi)} \hat{f}(\xi) \d \xi,
\end{equation}
and with phases obeying \eqref{eq:nrphase}.

In our numerical example, we consider the phase $\Phi_+$ and choose
\begin{align*}
r_1(x) & = \frac{1}{9}(2+\sin(4\pi x_1))(2+\sin(4\pi x_2)),\\
r_2(x) & = \frac{1}{9}(2+\cos(4\pi x_1))(2+\cos(4\pi x_2)).
\end{align*}
In each wedge $W_\ell$, the phase is then linearized and a low rank
separated approximation $UT$ of the matrix 
\[
A = \left( e^{2\pi i \Phi_\ell(x,\xi)} \right)_{x\in X,\xi\in W_\ell}
\]
is computed.  To estimate the approximation error, we randomly select
two sets $\Gamma$ and $\Delta$ of $s$ rows and $s$ columns. Put
$A_{\Gamma\Delta}$ to be the $s$ by $s$ the submatrix of $A$ with
these rows and columns. The separated rank approximation to
$A_{\Gamma\Delta}$ is then obtained by multiplying $U_\Gamma$ and
$T_\Delta$ where $U_\Gamma$ is the submatrix of $U$ with rows in
$\Gamma$ and $T_\Delta$ is that of $T$ with columns in $\Delta$. The
error is then estimated via 
\[
\frac{ \|A_{\Gamma\Delta} - U_\Gamma T_\Delta \|_F }
{ \|A_{\Gamma\Delta}\|_F }, 
\]
where $\|\cdot\|_F$ stands for the Frobenius norm. In our numerical
test, we set $s$ to be 200, and Table \ref{tbl:checkprob} displays
approximation errors for different combinations of problem size $N$
and accuracy $\eps$. The results show that the randomized approach
works quite well and that the estimated error is controlled well below
the threshold $\eps$.
\begin{table}[ht!]
  \begin{center}
        \begin{tabular}{|c|cccc|}
          \hline
              & $\eps=$1e-3 & $\eps=$1e-4 & $\eps=$1e-5 & $\eps=$1e-6 \\
          \hline
          $N=64$  & 3.57e-04 & 4.93e-05 & 3.21e-06 & 5.17e-07 \\
          $N=128$ & 3.11e-04 & 2.28e-05 & 4.19e-06 & 5.81e-07 \\
          $N=256$ & 2.85e-04 & 2.83e-05 & 2.94e-06 & 4.13e-07 \\
          $N=512$ & 1.66e-04 & 2.82e-05 & 4.38e-06 & 6.80e-07 \\
          \hline
        \end{tabular}
  \end{center}
  \caption{Relative errors of the low rank separated representation
        constructed using the randomized approach.}
  \label{tbl:checkprob}
\end{table}

Next, consider the relationship between the separation rank and the
threshold $\eps$. Corollary \ref{teo:lexing} shows that $\eps$ scales
like $N^{-p}$ for a fixed constant $p$ provided that the separation
rank grows gently like $p \log N$. In this experiment, we use the same
phase function $\Phi(x,\xi)$ in \eqref{eq:nrphase}, and show the
separation rank for different values of $N$ and $p$ in Table
\ref{tbl:checkrank}. These results suggest that the separation rank is
roughly proportional to both $p$ and the logarithm of $N$, which is
compatible with the theoretical estimate. Moreover, when $N$ is fixed,
the rank seems to grow linearly with respect to $p$, which possibly
implies that the constant $C(p)$ in Theorem \ref{teo:lexing} in fact
grows linearly with respect to $p$.

\begin{table}[ht!]
  \begin{center}
        \begin{tabular}{|c|ccccc|}
          \hline 
              & $p=$1 & $p=$1.5 & $p=$2 & $p=$2.5 & $p=$3 \\
          \hline
          $N=$64  & 7    &10    &14    &18    &22\\
          $N=$128 & 9    &12    &17    &21    &24\\
          $N=$256 & 9    &12    &17    &21    &25\\
          $N=$512 & 10   &15    &19    &24    &27\\
          \hline
        \end{tabular}
  \end{center}
  \caption{Ranks of the separated representation generated by the
    randomized approach for different values of $N$ and $p$.
    The prescribed error is equal to $N^{-p}$.
  }
  \label{tbl:checkrank}
\end{table}

We now turn to the numerical evaluation of $(L f)(x)$, 
\begin{equation}
  (Lf)(x) = \frac{1}{N} \sum_{\xi} e^{2\pi i \Phi(x,\xi)} \hat{f}(\xi), 
  \label{eq:dfioex}
\end{equation}
where the phase function $\Phi$ is the same as in
\eqref{eq:nrphase}. In this example, $f$ is an array of independently
and identically mean-zero normal random variables (Gaussian white
noise), which in some ways is the most challenging input. The
threshold $\eps$ is set to be $10 \, N^{-2}$ (i.e., $p=2$). To
estimate the error, we first pick $s$ points $\{x_i: i=1,\ldots,s\}$
from $X$ and put $\{ \widetilde{(Lf)}(x_i) \}$ for the output of our
algorithm (Section \ref{sec:algo-eval}). We then compare the values of
$\widetilde{(Lf)}(x_i)$ at these points with those of $\{ (Lf)(x_i)\}$
obtained by evaluating \eqref{eq:dfioex} directly. Finally, we
estimate the relative error with
\[
\sqrt{
  \frac{ \sum_i |(Lf)(x_i) - \widetilde{(Lf)}(x_i) |^2 }
  { \sum_i |(Lf)(x_i) |^2 }
}.
\]
Here, we choose $s=100$, and Table \ref{tbl:result_a} summarizes our
findings for various values of $N$.  The results show that our
algorithm performs well. The error is controlled well below threshold
and the speedup over the naive algorithm is significant for large
values of $N$.
\begin{table}[ht!]
  \begin{center}
        \begin{tabular}{|c|ccccc|}
          \hline
          $(N,\eps)$ & Preprocessing(s) & Evaluation(s) & Speedup & Error & Storage(MB) \\
          \hline
          (64,2.44e-03) &2.06e+00 &3.89e+00 &2.05e+00 &2.08e-03 &0.76\\
          (128,6.10e-04) &1.09e+01 &2.45e+01 &6.58e+00 &8.02e-04        &1.26\\
          (256,1.53e-04) &8.10e+01 &1.65e+02 &1.67e+01 &1.00e-04        &2.01\\
          (512,3.81e-05) &4.67e+02 &9.88e+02 &4.46e+01 &4.22e-05        &3.06\\
          \hline
        \end{tabular}
  \end{center}
  \caption{Numerical evaluation of $Lf(x)$ with $f$ a 
two dimensional white-noise array. The second and third
    columns give the number of seconds spent in the preprocessing and
    evaluation steps respectively. The fourth column shows the speedup
    factor over the naive algorithm for computing $(Lf)(x)$ using 
    the direct summation  \eqref{eq:dfioex}. The fifth column is the 
    estimated relative error and the
    last gives the amount of memory used in terms of megabytes.}
  \label{tbl:result_a}
\end{table}

We have only considered the evaluation of FIOs in ``Egorov'' form thus
far (constant amplitude) but the algorithm described in Section
\ref{sec:algo} can be easily extended to operate with general
amplitudes provided that the term $a(x,\xi)$ also admits a low rank
separated representation in the variables $x$ and $\xi$.

To study the performance of our algorithm in the more general setup of
variable amplitudes, we continue with the example where $f$ is
integrated along ellipses \eqref{eq:ellipses_FIO} (recall the phase
\eqref{eq:nrphase}). The Appendix shows that a possible choice for the
amplitudes $a_\pm(x,\xi)$ and phases $\Phi_\pm(x,\xi)$ is
\begin{align}
  \label{eq:ellipses_amp}
  a_\pm(x,\xi) & = \frac{1}{4 \pi} \left( \, J_0(2 \pi \rho(x,\xi)) \pm i
    Y_0(2 \pi \rho(x,\xi)) \, \right) e^{\mp 2 \pi i \rho(x,\xi)},\\ 
  \Phi_\pm(x,\xi) & = x\cdot\xi \pm \rho(x,\xi)
\end{align}
with
\[
\rho(x,\xi)  = \sqrt{r_1^2(x) \xi_1^2 + r_2^2(x) \xi_2^2}.
\]
Here, $J_0$ and $Y_0$ are Bessel functions of the first and second
kind respectively, see the Appendix for details. 

For the axes lengths, set
\begin{equation}
  r_1(x) = r_2(x) \equiv r(x) = 
\frac{1}{16} (3+\sin(4\pi x_1))(3+\sin(4\pi x_2))
\label{eq:circlerad}
\end{equation}
(which means that our ellipses are circles). We compute $(L_+ f)(x)$
for different values of $N$ and $\eps$ and provide the results in
Table \ref{tbl:result_b}. The computational analysis shows that our
algorithm performs equally well in the variable amplitude case.  For
$N = 512$, the speedup factor over the naive evaluation is about 162.
\begin{table}[ht!]
  \begin{center}
        \begin{tabular}{|c|ccccc|}
          \hline
          $(N,\eps)$ & Preprocessing(s) & Evaluation(s) & Speedup & Error & Storage(MB) \\
          \hline
          (64,2.44e-03)  &2.18e+01 &3.67e+01 &4.54e+00 &7.30e-04 &0.37 \\
          (128,6.10e-04) &1.09e+02 &1.65e+02 &1.49e+01 &4.00e-04 &0.59 \\
          (256,1.53e-04) &6.62e+02 &8.46e+02 &4.49e+01 &1.39e-04 &0.89 \\
          (512,3.81e-05) &3.42e+03 &4.43e+03 &1.62e+02 &3.69e-05 &1.38 \\
          \hline
        \end{tabular}
  \end{center}
  \caption{Numerical evaluation of $Lf(x)$ with $f$ a 
two dimensional white-noise array. }
  \label{tbl:result_b}
\end{table}

\begin{figure}
  \begin{center}
        \begin{tabular}{ccc}
          \includegraphics[width=1.9in]{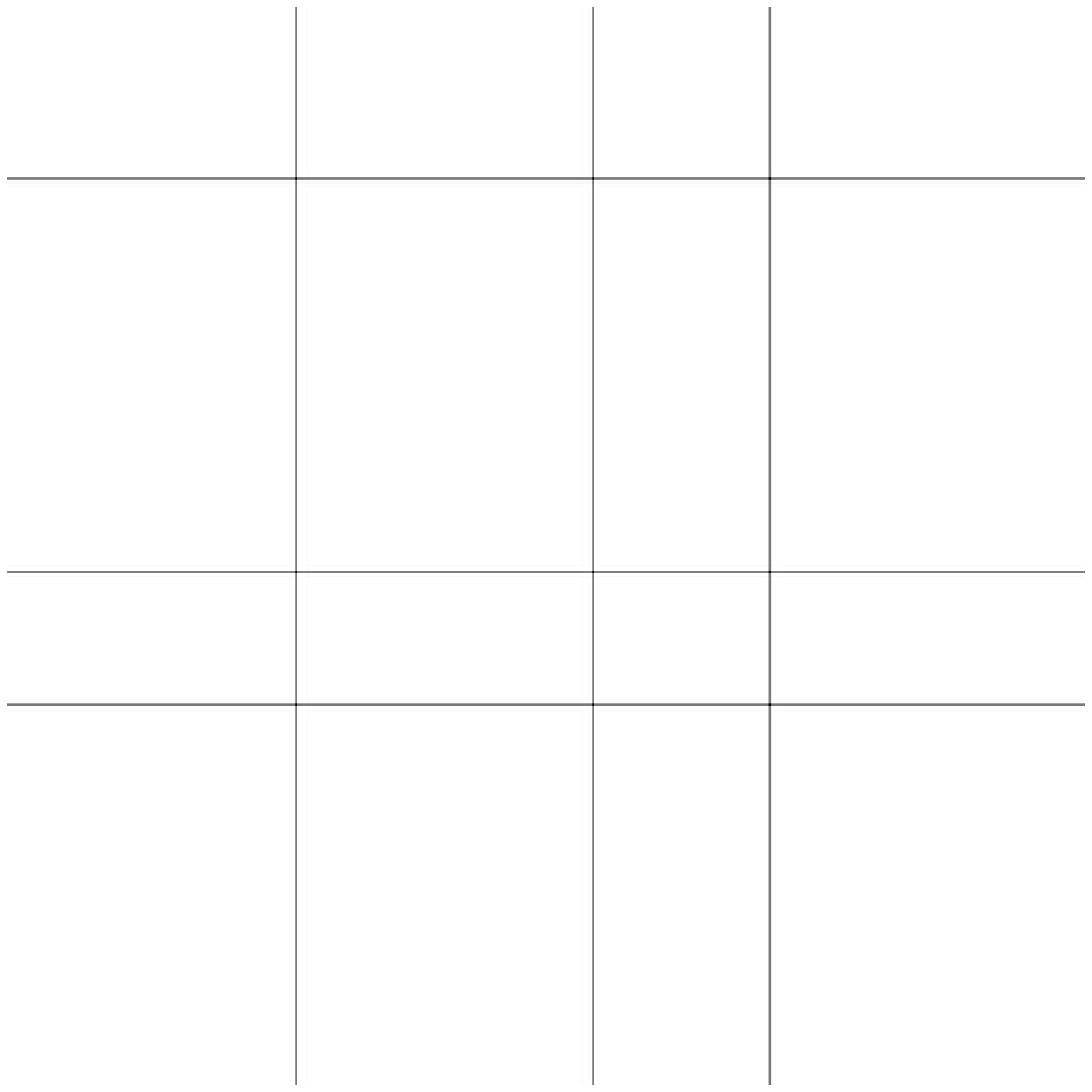}&
          \includegraphics[width=1.9in]{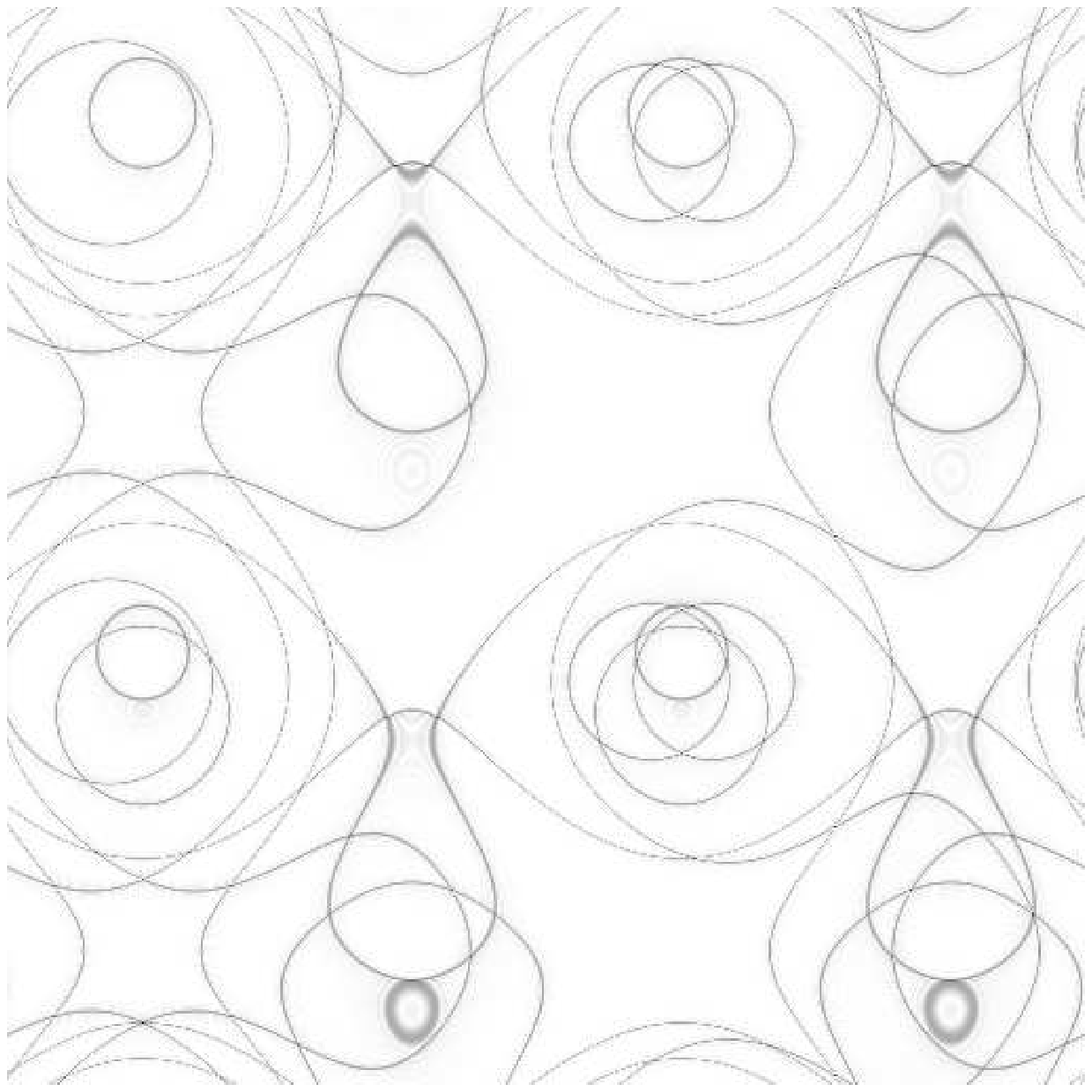}&
          \includegraphics[width=1.9in]{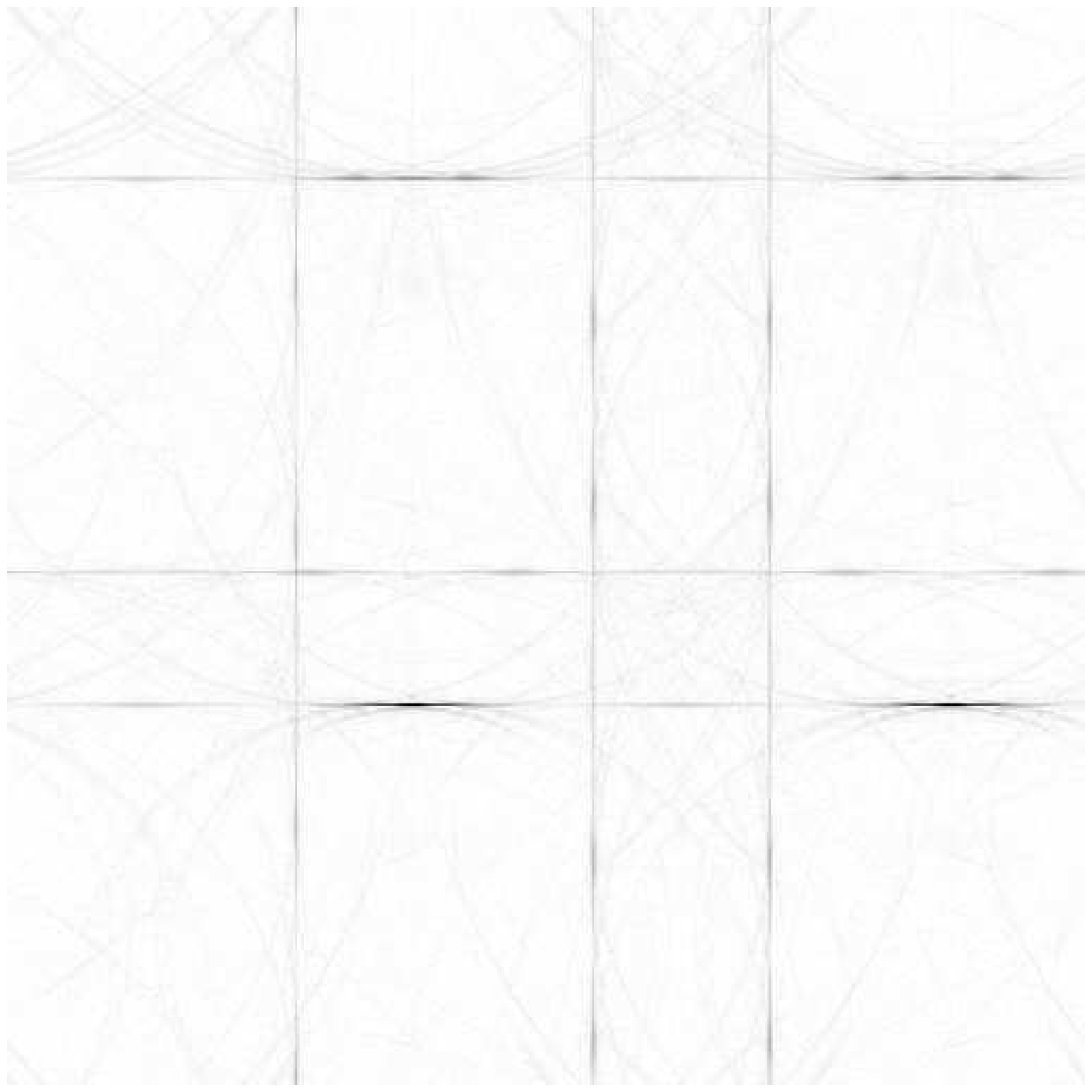}\\
          \includegraphics[width=1.9in]{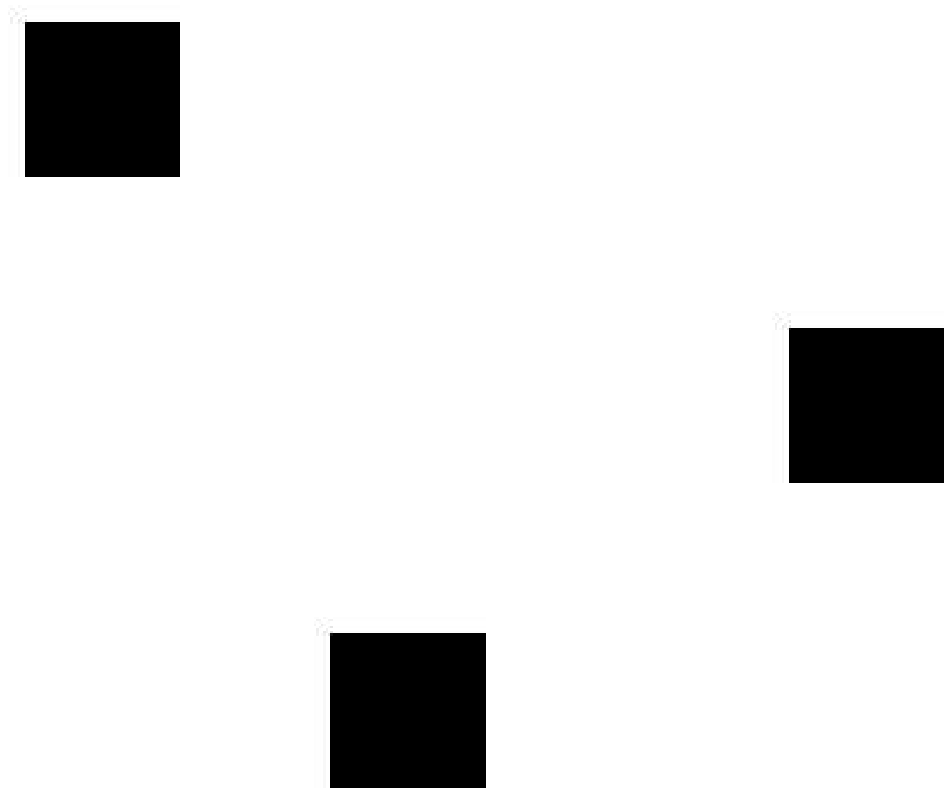}&
          \includegraphics[width=1.9in]{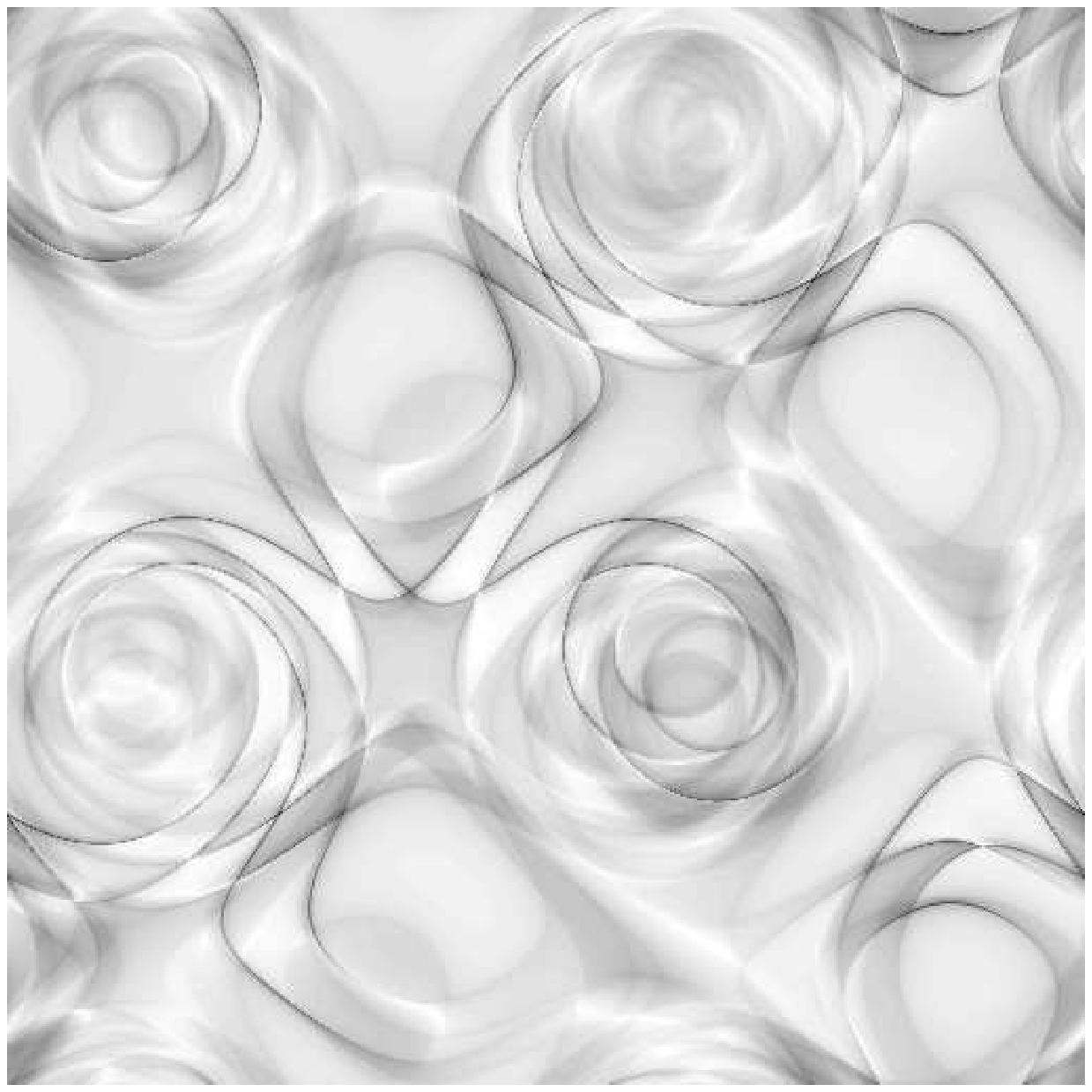}&
          \includegraphics[width=1.9in]{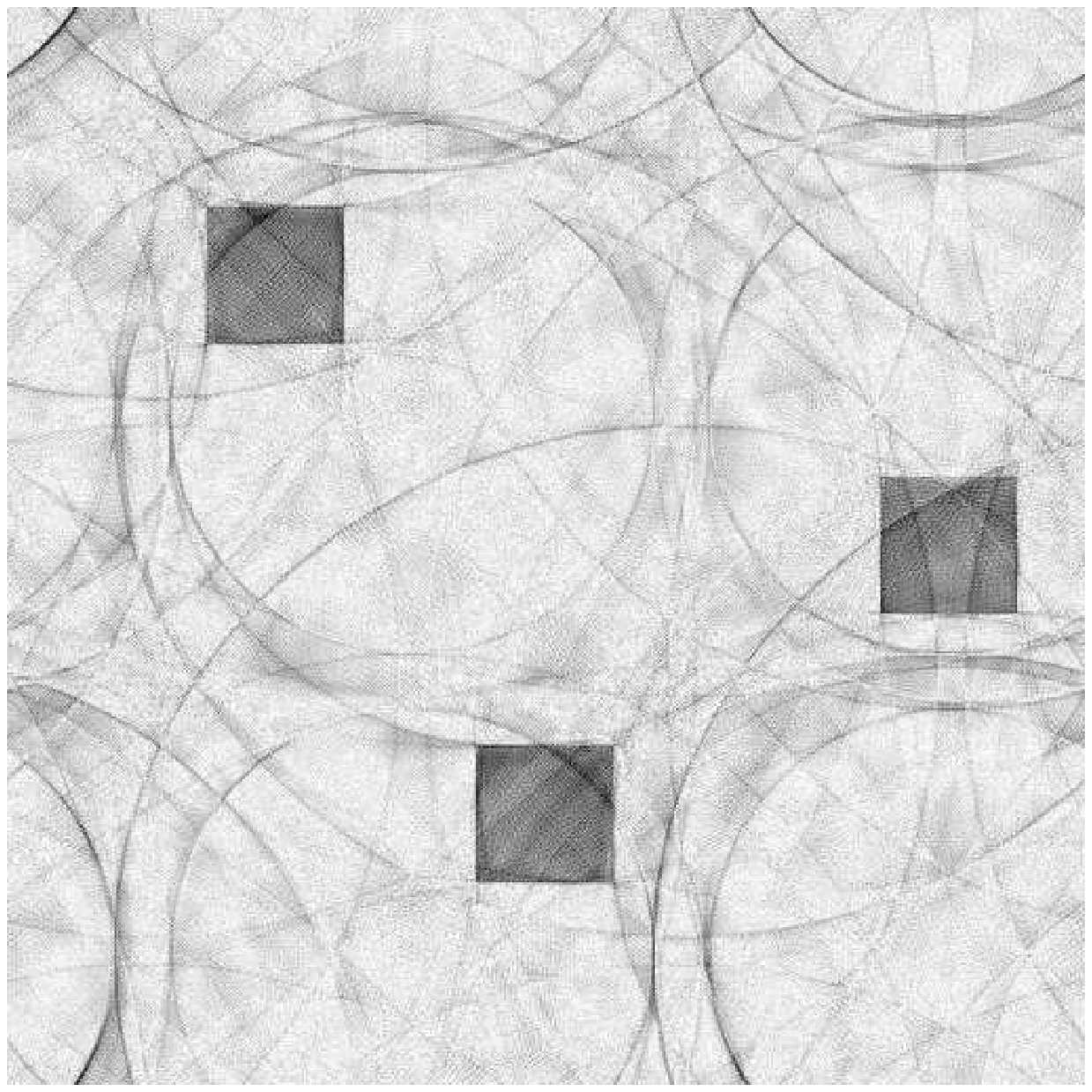}\\
          \includegraphics[width=1.9in]{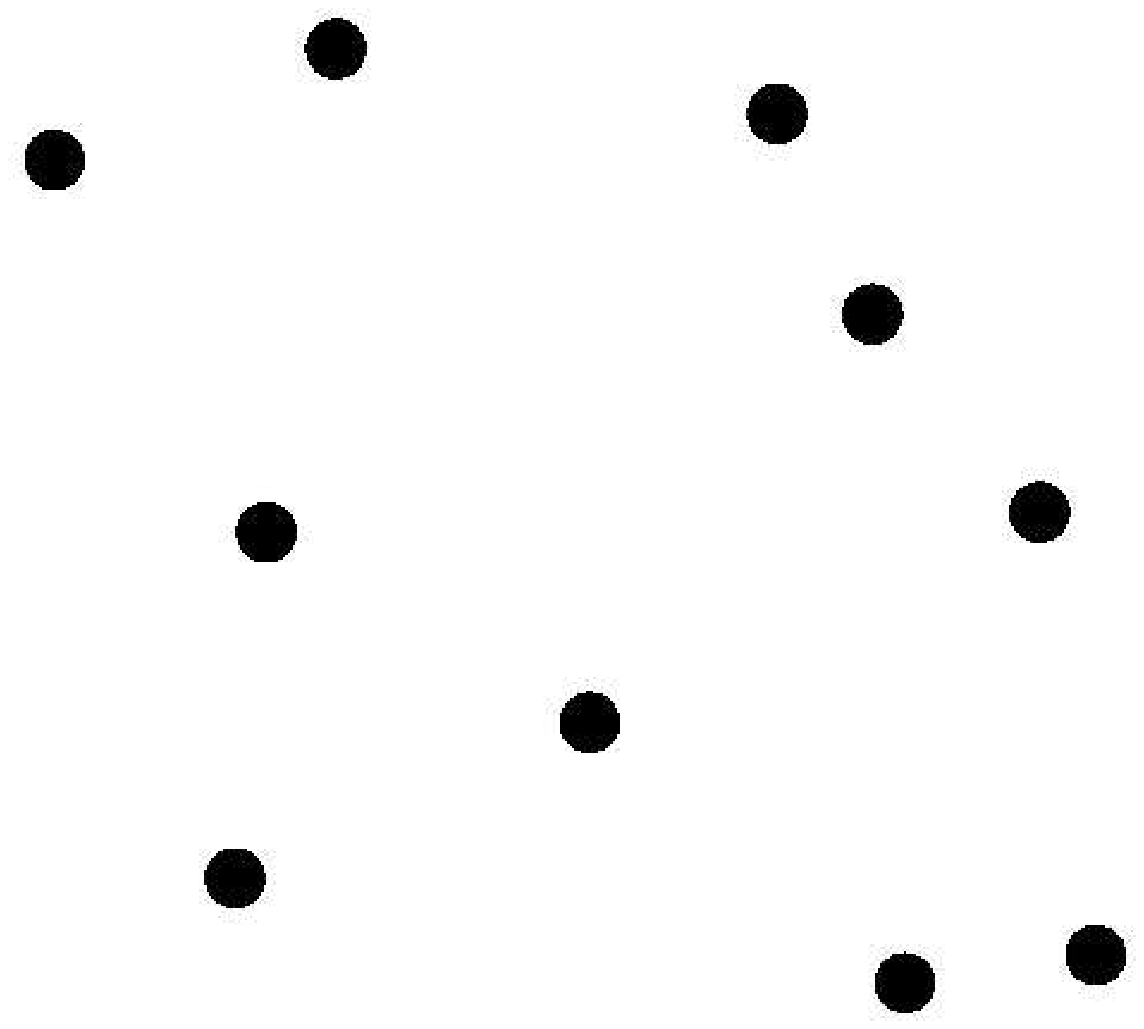}&
          \includegraphics[width=1.9in]{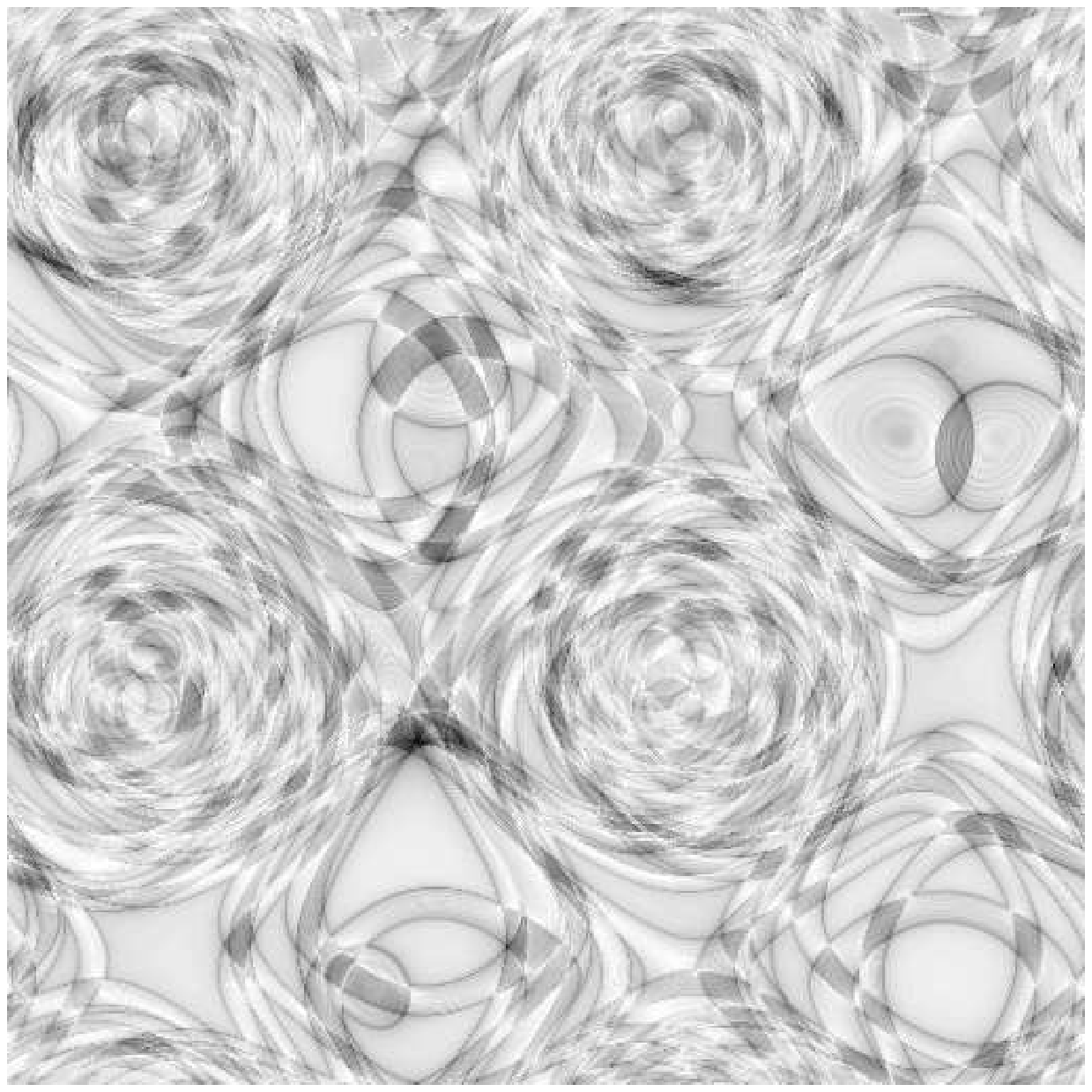}&
          \includegraphics[width=1.9in]{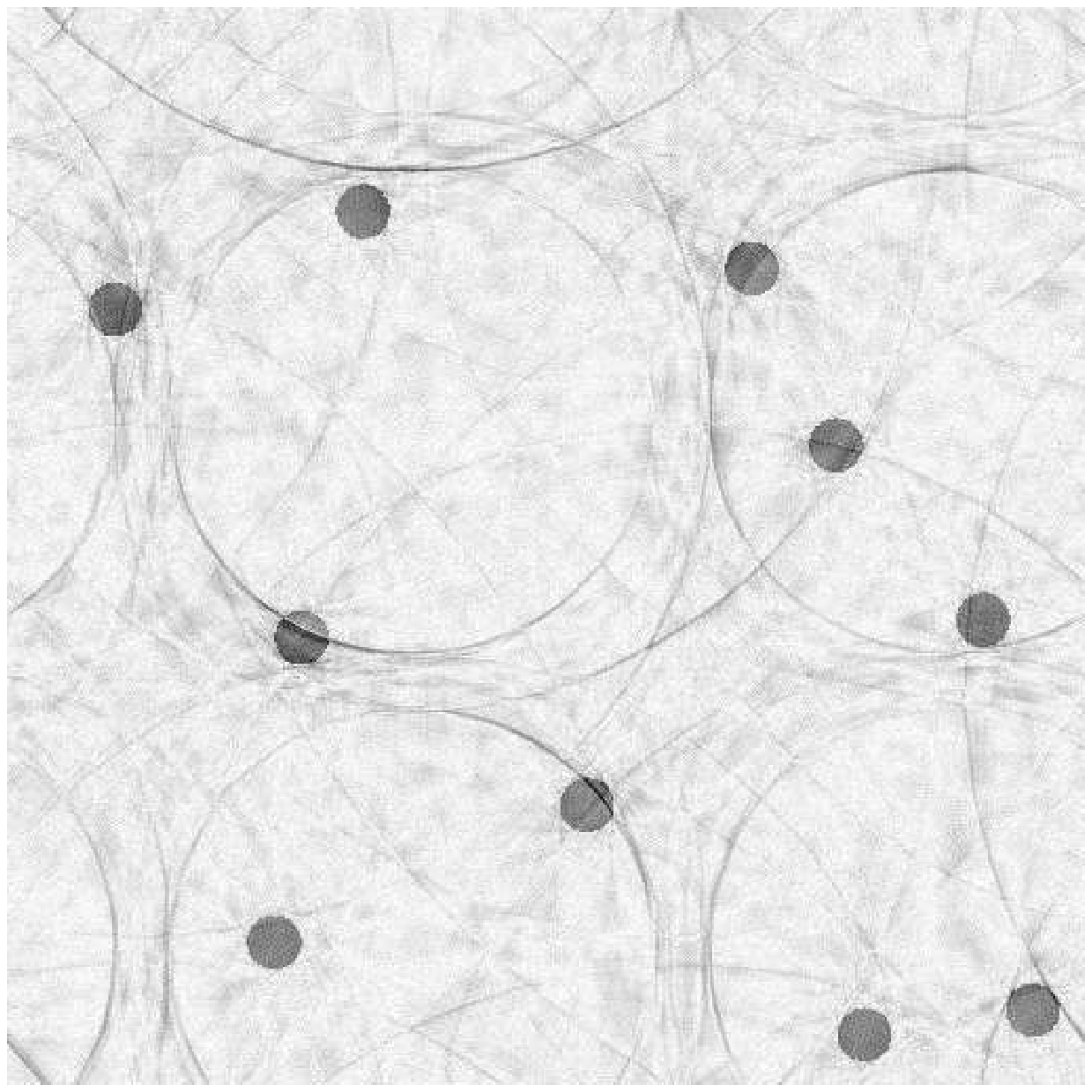}
        \end{tabular}
  \end{center}
  \caption{Numerical verification of the fact $WF(L^*L f) =
    WF(f)$. Each row, from left to right, shows the magnitudes of
    $f(x)$, $(Lf)(x)$ and $(L^* L f)(x)$ .  Notice that the wave front
    set of $f(x)$ and $(L^* L f)(x)$ are numerically as close as they
    can be. \emph{Remark}: the images on the left column and on the
      right column are not supposed to be the same; only their
      ``singularities'' coincide. In other words, the adjoint $L^*$ is
      not the inverse of $L$.}
  \label{fig:wf}
\end{figure}

An extremely important property of Fourier integral operators is that,
under the nondegeneracy condition
\[
\mbox{det} \left( \frac{\pd^2 \Phi}{\pd x_i \pd x_j} \right) \ne 0,
\]
the composition of an FIO with its adjoint preserves the singularities
of the input function. Mathematically speaking, if $WF(f)$ is the wave
front set of $f$ \cite{duistermaat-book,symes-1998-mrs}, then
\[
WF(L^* L f) = WF(f). 
\]
This property serves as the foundation for most of the current imaging
techniques in reflection seismology \cite{symes-1998-mrs}. In the
final example of this section, we verify this phenomenon
numerically. We choose the phase function to be
\[
\Phi(x,\xi) = x\cdot \xi + r(x)|\xi|, 
\]
where $r(x)$ is given by \eqref{eq:circlerad}, and compute $(L^* L f)$
using the algorithm discussed in Sections \ref{sec:algo-eval} and
\ref{sec:algo-evaladj}. Figure \ref{fig:wf} displays results for three
input functions with different kinds of singularities. Looking at the
picture, we see that the singularities of $L f$ are of course
different than those of $f$, but we also see that the singularities of
$L^* L f$ coincide with those of $f$.

\section{Discussion}
\label{sec:disc}

\subsection{About randomized algorithms}

The method used in the randomized preprocessing step was first
introduced by Kapur and Long \cite{kapur-1997-IES3}. Lately, there has
been a lot of research devoted to the development of randomized
algorithms for generating low rank factorizations, and we would like
to discuss some of this work.

Drineas, Kannan and Mahoney \cite{drineas-2006-fmcam2} describe a
randomized algorithm for computing a low-rank approximation to a fixed
matrix. The main idea is to form a submatrix by selecting columns with
a probability proportional to their norm. Since this work is about
unstructured general matrices, it does not guarantee a small
approximation error. As an example, suppose all the columns of the
matrix have the same norm and one of them is orthogonal to the span of
the other columns. Unless this column is selected, the orthogonal
component is lost and the resulting approximation is poor.

Our situation is different. Since each entry of our matrix
\[
A = \left( e^{2\pi i \Phi_\ell(x,\xi)} \right)_{x\in X,\xi\in W_\ell}
\]
has unitary magnitude, the uniform probability used in our algorithm
is actually the same as that proposed above
\cite{drineas-2006-fmcam2}. In some ways then, our approach is a
special case of that of Drineas et.~al. But the point is that our
matrix has a special structure. As we argued earlier, the columns of
$A$ are often highly correlated and we believe that this is the reason
why the randomized subsampling performs well.

A recent article by Martinsson, Rokhlin and Tygert
\cite{martinsson-2006-rafam} presents a new randomized solution to the
same problem. 
The only inconvenience of this algorithm, probably inevitable for
general matrices, is that one needs to visit all the entries of the
matrix multiple times. This can be quite costly in our setup since
there are $O(N^4)$ entries. This is why we adopt the method by Kapur
and Long.





\subsection{Storage compression}

We would like to comment on the storage compression strategy discussed
at the end of Section \ref{sec:algo-rand}. In fact, what we described
there can be viewed as a new way of compressing low rank matrices.

In a general context, the entries of a matrix can be viewed as
interaction coefficients between a set of objects indexed by the rows
and another set indexed by the columns. In our case, the first set
contains the grid points $x$ in $X$, while the second set consists of
the frequencies $\xi$ in $W_\ell$. Call these two sets $I$ and
$J$, and the interaction matrix $A_{I,J}$.  The standard practice for 
compressing $A_{I,J}$ is to find two sets $I'$ and $J'$ of
smaller sizes and form an approximation
\[
A_{I,J} \approx M_{I, I'} M_{I', J'} M_{J', J}.
\]
Here $I'$ is either a subset of $I$ or a set which is close by in some
sense, and likewise for $J'$ and $J$. For example, in the fast
multipole method of Greengard and Rokhlin \cite{greengard-1987-FAPS},
$J'$ is the multipole representation at the center of the box
containing $J$ while $I'$ is the local representation at the center of
the box containing $I$. The matrices $M_{I, I'}$, $M_{I', J'}$ and
$M_{J', J}$ are implemented as the multipole-to-multipole,
multipole-to-local and local-to-local translations. This becomes even
more obvious when one considers the newly proposed kernel independent
fast multipole method by Ying, Biros and Zorin
\cite{ying-2004-KIAFMA}. There, $I'$ and $J'$ are the equivalent
densities supported on the boxes containing $I$ and $J$, while $M_{I,
  I'}$, $M_{I', J'}$ and $M_{J', J}$ can be computed directly from
interaction matrices and their inverses. In both cases, we are
fortunate in the sense that prior knowledge offers us efficient ways
to multiply $M_{I, I'}$, $M_{I', J'}$ and $M_{J', J}$ with arbitrary
vectors. Whenever this is not true, one might be forced to store these
matrices, which could be quite costly.

\begin{figure}
  \begin{center}
    \includegraphics[height=2in]{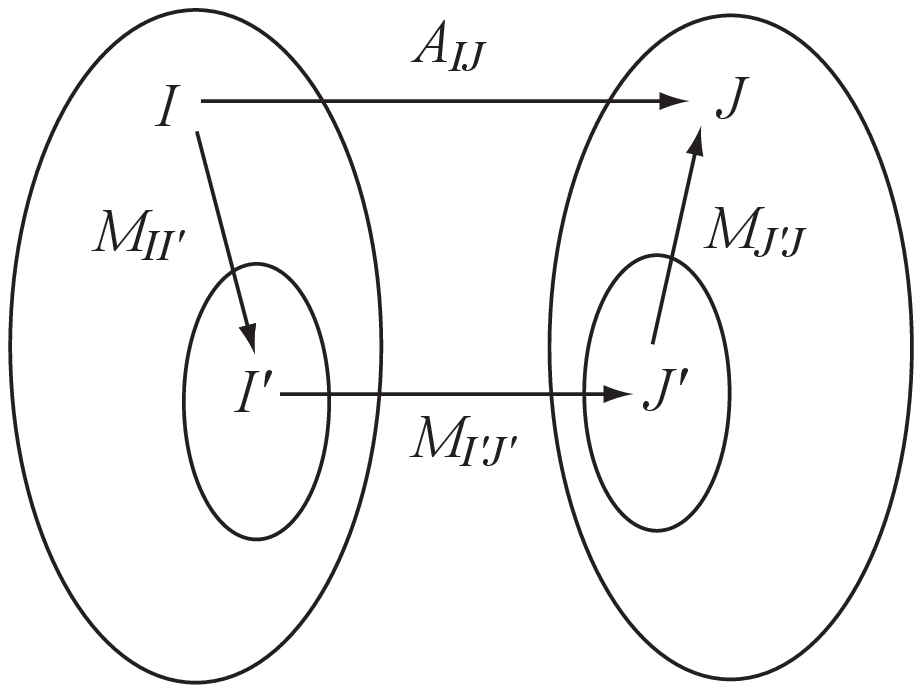}
    \includegraphics[height=2in]{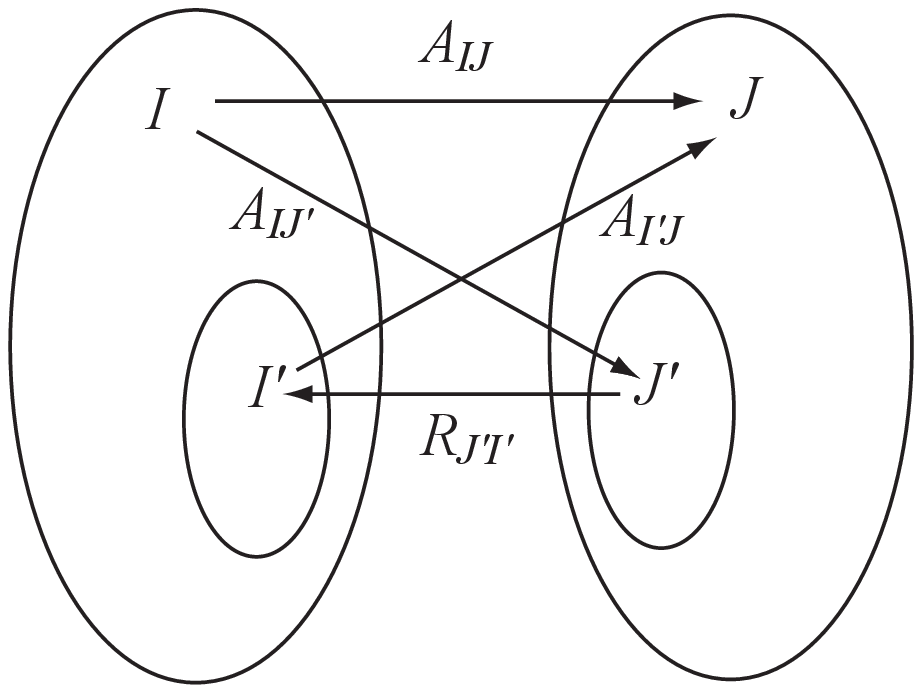}
  \end{center}
  \caption{Factorization of interaction between $A$ and $B$. (a) the
    standard scheme, (b) the scheme abstracted from the storage
    compression method used \eqref{eq:probnewfac}.}
  \label{fig:storage}
\end{figure}

What we have presented in \eqref{eq:probnewfac} is a totally different
factorization:
\[
A_{IJ} \approx A_{I J'} R_{J' I'} A_{I' J}.
\]
Notice that since $A_{I J'}$ and $A_{I' J}$ are interaction matrices
themselves, there is no need to store them as long as we can compute
the interaction coefficients easily. The only thing we need to keep in
storage is the matrix $R_{J' I'}$. However, as long as the interaction
is low rank, $I'$ and $J'$ have far fewer objects than $I$ and $J$, so
that $R_{J' I'}$ only uses very little storage. Finally, we would like
to point out that, instead of representing the interaction from $J'$
(a subset of $J$) to $I'$ (a subset of $I$), $R_{J' I'}$ is a reverse
interaction. Figure \ref{fig:storage} shows conceptually how the new
factorization differs from the standard one.

\subsection{Curvelets, wave atoms and beamlets}

There might be other ways of evaluating Fourier integral operators,
and we would like to discuss their relationships with the approach
taken in this paper.

Curvelets, proposed by Cand\`{e}s and Donoho \cite{candes-2004-ntf},
are two dimensional waveforms which are highly anisotropic in the fine
scales.  Each curvelet is identified with three numbers to indicate
its scale, orientation and position, and the set of all curvelets form
a tight frame. Recently, Cand\`{e}s and Demanet
\cite{candes-2003-cfio,candes-2005-crwpos} have shown that the
curvelet representation of the Fourier integral operators is optimally
sparse. More precisely, a Fourier integral operator only has $O(N^2)$
nonnegligible entries in the curvelet domain. The wave atom frame,
which is recently introduced by Demanet and Ying
\cite{demanet-2006-wasop}, has the same property. If we were able to
find such a representation efficiently, we would hold an $O(N^2\log
N)$ algorithm for evaluating a Fourier integral operator which would
operate as follows:
\begin{enumerate}
\item Apply the forward curvelet transform to the input and get
  curvelet coefficients.
\item Apply the sparse FIO to the curvelet coefficient sequence. 
\item Apply the inverse curvelet transform.
\end{enumerate}
Both steps 1 and 3 require at most $O(N^2\log N)$ operations
\cite{FDCT}.

Constructing the curvelet representation of FIO from the phase
function $\Phi(x,\xi)$ efficiently has, however, proved to be
nontrivial. At the moment, we are only able to construct an
approximation which is asymptotically accurate by studying the
canonical relation embedded inside the phase function
$\Phi(x,\xi)$. Such a construction would be adequate if we were
interested in applying an FIO to input functions with only high
frequency modes. However, one often wants a representation which is
accurate for all frequency modes, and we are currently not aware of
any efficient method for constructing such a representation.

Beamlets \cite{donoho-2002-bmia} were introduced by Donoho and Huo at
roughly the same time as curvelets. Beamlets are small segments at
different positions, scales and orientations. As pointed out in
Section \ref{sec:results}, curvilinear integrals make up an important
subclass of FIOs, and beamlets may offer ways to efficiently compute
such simpler integrals. One might think of something like this:
\begin{enumerate}
\item Compute the beamlet coefficient sequence of the input.
\item For each $x\in X$, figure out the integration curve and
  approximate it with a chain of beamlet segments. Sum up the beamlet
  coefficients along the chain.
\end{enumerate}
Assuming the integration curves are twice differentiable, we would
need about $\sqrt{N}$ beamlet segments to approximate each
curve. Thus, the overall complexity of this algorithm might scale like
$O(N^{2.5})$, which is the same scaling as that of our algorithm. The
problem is that it is unclear how one would efficiently approximate
the integration curve with beamlet segments without sacrificing
accuracy. Situations in which the input function $f$ is highly
oscillatory or in which the integration curves have parts with a high
curvature seem very problematic.

Our algorithms decompose the FIO in the frequency domain whereas the
beamlet based approach processes data in the spatial
domain. Sandwiched right in the middle, curvelets and wave atoms
operate in the phase-space---the product of the frequency and of the
spatial domains. We believe that operating in phase-space by
exploiting the microlocal properties of FIOs would be important to
bring down the complexity to the optimal value of about $N^2$
operations.

\appendix
\section{Integration Along Ellipses}

The material in this section is probably not new, but we expand on it
for the convenience of the nonspecialist. Consider the generalization
Radon transform that consists in integrating $f(x)$ along ellipses of
axes lengths $r_1(x)$ and $r_2(x)$, and centered around $x$:
\[
Gf(x) = \int f \left( x + \begin{pmatrix} r_1(x) \cos \theta \\ r_2(x) \sin \theta
\end{pmatrix} \right) \, d\theta.
\]
We want to recast it as a sum of FIOs. Let us start by writing
\[
Gf(x) = \int K(x,\xi) \hat{f}(\xi) \, d\xi,
\]
with
\[
K(x,\xi) = e^{2 \pi i x \cdot \xi} \int \exp \left[ 2 \pi i
  \begin{pmatrix} r_1(x) \cos \theta \\ r_2(x) \sin \theta
  \end{pmatrix} \cdot \xi \right] \, d\theta.
\]
Put $\rho(x,\xi) = \sqrt{r_1^2(x) \xi_1^2 + r_2^2(x) \xi_2^2}$ and rewrite
\[
K(x,\xi) = e^{2 \pi i x \cdot \xi} \int \exp\left[ 2 \pi i \rho(x,\xi)
  \begin{pmatrix} \cos \theta \\ \sin \theta \end{pmatrix} \cdot  \begin{pmatrix} \alpha(x,\xi) \\ \beta(x,\xi) \end{pmatrix} \right] \, d\theta.
\]
Here $\alpha^2 + \beta^2 = 1$, and both $\alpha$ and $\beta$ depend on
$x$ and $\xi$ but the value of the integral is independent of their
particular value. This is because any change of variables $\theta \to
\theta + \phi(x,\xi)$, effectively corresponding to a rotation of the
unit vector $(\alpha, \beta)$, keeps the integral invariant. So we
may as well take $\alpha = 1$, $\beta = 0$ and obtain
\[
K(x,\xi) = e^{2 \pi i x \cdot \xi} \int e^{2 \pi i \rho(x,\xi) \cos \theta} \, d\theta = \frac{e^{2 \pi i x \cdot \xi}}{2 \pi} J_0(2 \pi \rho(x,\xi)).
\]
Of course the Bessel function $J_0$ oscillates, and we need to extract
the phase from its asymptotic behavior
\[
J_0(2 \pi \rho(x,\xi)) \sim \sqrt{\frac{1}{\pi^2 \rho(x,\xi)}} \cos \left( 2 \pi \rho(x,\xi) -
\frac{\pi}{4} \right).
\]
The idea is now to express $J_0(2\pi \rho(x,\xi))$ as a sum of
two terms, each of which being the product between a smooth amplitude
(a demodulated version of $J_0$ or the envelope of $J_0$ if you will)
and the oscillatory exponential $e^{\pm 2 \pi i \rho(x,\xi)}$.  In
effect, $K$ is decomposed as a sum of two FIOs:
\[
K(x,\xi) = a_+(x,\xi) e^{2 \pi i \Phi_+ (x,\xi)} + a_-(x,\xi) e^{2 \pi i \Phi_-
  (x,\xi)},
\]
with
\[
\Phi_\pm(x,\xi) = x \cdot \xi \pm \rho(x,\xi).
\]
There are different ways to choose the amplitudes. One way is to let
$Y_0$ be the Bessel function of the second kind of order zero
\cite{abramowitz-1972-hmf} and exploit the identity $2 J_0 = (J_0 + i
Y_0) + (J_0 - i Y_0)$, which allows to write
\[
a_\pm(x,\xi) = \frac{1}{4 \pi} \left( \, J_0(2 \pi \rho(x,\xi)) \pm i
Y_0(2 \pi \rho(x,\xi)) \, \right) e^{\mp 2 \pi i \rho(x,\xi)}.
\]

Both amplitudes behave asymptotically like $\sqrt{1/\pi^2
  \rho(x,\xi)}$ as $x \to \infty$, which incidentally shows that the
order of the FIO is $-1/2$. The logarithmic singularity of $Y_0$ near
the origin in $\xi$ is mild and easily regularized with no loss of
accuracy.

\end{document}